\documentclass{amsart}

\pdfoutput=1
\usepackage[utf8]{inputenc}
\usepackage[T1]{fontenc}
\usepackage{lmodern}
\usepackage{mathtools}
\allowdisplaybreaks
\usepackage{amssymb}
\usepackage{amsthm}
\usepackage[usenames,dvipsnames,svgnames,table]{xcolor}
\usepackage{embrac}
\ChangeEmph{(}[-.01em,.04em]{)}[.04em,-.05em]
\usepackage{textcase}
\usepackage{enumitem}
\usepackage[acronym]{glossaries}
\usepackage{graphicx}
\usepackage{float}
\usepackage{booktabs}
\usepackage{csquotes}
\usepackage{multirow}
\usepackage{tikz}
\usepackage{subcaption}
\tikzset{node distance=2cm, auto}
\usetikzlibrary{angles,quotes,calc,positioning,shapes,matrix,arrows}
\usepackage{tikz-cd}
\usetikzlibrary{}
\usepackage{tkz-euclide}
\usepackage{etoolbox}
\usepackage{url}
\usepackage[hyperindex=true, colorlinks=true, linkcolor=Maroon, citecolor=DarkGreen, urlcolor=NavyBlue, breaklinks=true,linktocpage=true,linktoc=all]{hyperref}
\usepackage[english,capitalise,noabbrev]{cleveref}
\crefdefaultlabelformat{#2\textup{#1}#3}
\crefname{equation}{equation}{equations}
\usepackage[logical-quotes,backrefs,initials]{amsrefs}

\newcommand{\R}{\mathbb{R}}

\renewcommand{\S}{\mathbb{S}}

\newcommand{\RP}{\mathbb{RP}}

\newcommand{\E}{\mathcal{E}}

\newcommand{\greek}[1]{\mathrm{#1}}

\newcommand{\Beta}{\greek{B}}

\newcommand{\contained}{\subset}
\newcommand{\containedeq}{\subseteq}

\newcommand{\union}{\cup}

\newcommand{\suchthat}{\, : \,}
\newcommand{\st}{\suchthat}

\newcommand{\defined}{\coloneqq}

\newcommand{\SO}{\operatorname{SO}}

\newcommand{\Id}{\operatorname{Id}}
\newcommand{\diag}{\operatorname{diag}}
\newcommand{\trace}{\operatorname{tr}}

\newcommand{\tendsto}{\rightarrow}

\newcommand{\from}{\colon}

\newcommand{\composed}{\circ}

\newcommand{\vol}{\operatorname{vol}}

\newcommand{\unif}{\operatorname{unif}}

\newcommand{\eval}[1][\bigg]{#1|}
\DeclarePairedDelimiterX\norm[1]\lVert\rVert{\ifblank{#1}{\:\cdot\:}{#1}}
\DeclarePairedDelimiterX\abs[1]\lvert\rvert{\ifblank{#1}{\:\cdot\:}{#1}}
\DeclarePairedDelimiterX\set[1]{\{}{\}}{\ifblank{#1}{\: \:}{#1}}
\DeclarePairedDelimiterX\innerprod[2]\langle\rangle
{\ifblank{#1#2}{\,\cdot,\cdot\,}
	{#1,#2}}
\DeclarePairedDelimiterX\floor[1]\lfloor\rfloor{\ifblank{#1}{\:\cdot\:}{#1}}
\DeclarePairedDelimiterXPP\expectation[2]{\ifblank{#1}{\mathbb{E}}{\mathbb{E}_{#1}}}\lparen\rparen{}{\ifblank{#2}{\:\cdot\:}{#2}}

\newcommand{\transpose}{^{\text{\tiny$\mathsf{T}$}}}

\newcommand{\Gegenbauer}[2]{C_{#1}^{(#2)}}

\newcommand{\Jacobi}[3]{P_{#1}^{(#2,#3)}}

\newcommand{\dd}{\mathop{}\!\mathrm{d}}

\newcommand{\BesselJ}[1]{J_{#1}}

\renewcommand{\digamma}{\psi}

\newcommand{\s}{\mathbb{S}^2}

\newcommand{\so}{\SO(3)}

\theoremstyle{plain}
\newtheorem{theorem}{Theorem}[section]

\newtheorem{lemma}[theorem]{Lemma}
\newtheorem{proposition}[theorem]{Proposition}

\theoremstyle{definition}

\newtheorem{definition}[theorem]{Definition}

\theoremstyle{remark}

\newtheorem{remark}[theorem]{Remark}

\setlist[enumerate,1]{label = \textup{(\textsc{\roman*})}}
\setlist[enumerate,2]{label = \textup{(\textit{\alph*})}}
\setlist[enumerate,3]{label = \textup{(\arabic*)}}

\subjclass[2020]{Primary 31C12, 60G55; Secondary 31B15, 52C35}
\keywords{Logarithmic Energy, Rotation group, Point arrangements, Random Polynomials}
  
\begin{document}
	
\title{Points on $\operatorname{SO}(3)$ with low logarithmic energy}

\author[C. Beltrán]{Carlos Beltrán} 
\address{Carlos Beltrán: Departamento de Matemáticas, Estadística y Computación, Universidad de Cantabria, Avda. Los Castros s/n, 39005 Santander, Spain}
\email{beltranc@unican.es}
\thanks{The first and fourth authors have been supported by grant PID2020-113887GB-I00 funded by MCIN/AEI/
	10.13039/501100011033. } 

\author[F. Carrasco]{Federico Carrasco}
\address{Federico Carrasco: Facultad de Ciencias Economicas y de Administración, Universidad de la República, Uruguay}
\email{federico.carrasco@fcea.edu.uy}

\author[D. Ferizović]{Damir Ferizović}
\address{Damir Ferizović: Department of Mathematics, KU Leuven, Celestijnenlaan 200b, Box 2400, 3001 Leuven, Belgium}
\email{damir.ferizovic@kuleuven.be}
\thanks{The third author is supported by the FWO grant 1231425N}

\author[P.\,R. López-Gómez]{Pedro R. López-Gómez}
\address{Pedro R. López-Gómez: Departamento de Matemáticas, Estadística y Computación, Universidad de Cantabria, Avda. Los Castros s/n, 39005 Santander, Spain}
\email{lopezpr@unican.es}
\thanks{The fourth author has also been supported by grant PRE2021-
	097772 funded by MCIN/AEI/ 10.13039/501100011033 and by “ESF Investing in your future”.} 
	
\date{\today{}}

\begin{abstract}
	We describe several randomized collections of $3\times 3$ rotation matrices and analyze their associated logarithmic energy. The best one (that is, the one attaining the lowest expected logarithmic energy) is constructed by choosing $r$ spherical points, which come from the zeros of a randomly chosen degree $r$ polynomial, and considering at each of these points a set of $s$ evenly distributed rotation matrices. This construction yields a new upper bound on the minimal logarithmic energy of $n=rs$ rotation matrices.
\end{abstract}
	
\maketitle

\section{Introduction and main result}

A central theme in approximation theory is the constructive characterization of point sets within a given space that exhibit specific, desirable properties. This process involves designing and analyzing point distributions that meet various criteria, such as minimizing interpolation error, maximizing density in certain regions, or ensuring uniform coverage of a domain. The most basic and famous examples are probably polynomial approximation, where Chebyshev nodes are used to minimize the Runge phenomenon, leading to more accurate interpolation, or numerical integration, where Gaussian quadrature selects points and weights to achieve exact integration for polynomials of a given degree. A lesser-known yet remarkable result is Fejér's solution to the so-called Fekete problem on the unit interval. The problem involves finding points $x_1,\dotsc,x_n\in[-1,1]$ that minimize the logarithmic energy, given by
\begin{equation}\label{eq:intervalo}
	\sum_{i\neq j}\log\abs{x_i-x_j}^{-1}.
\end{equation}
Fejér showed that the optimal configuration consists of the endpoints of the interval together with the $n-2$ zeros of a Gegenbauer polynomial. This solution highlights deep connections between potential theory and orthogonal polynomials; see \cite{Fejer,Granada} for further details.

In higher-dimensional spaces, the problem becomes significantly more complex, and far less is known about explicit, optimal point distributions. A substantial body of research has focused on the unit sphere, where optimality is often defined through the minimization of certain energies, in the spirit of \eqref{eq:intervalo}. In this context, the absolute value in the expression of the energy is replaced by the Euclidean distance, and various energy models are considered, such as the logarithmic energy, the sum of the inverses of the mutual distances, or, more generally, Riesz energies, where the distance is raised to an arbitrary power. Foundational work in this area includes \cite{Whyte1952}, and the literature has since grown extensively. While a comprehensive survey is beyond the scope of this paper, we mention some works such as \cite{RakhmanovSaffZhou1994,Brauchart2008,AlishahiZamani,BrauchartGrabner2015,BrauchartHardinSaff2012,HMS,BeltranMarzoOrtega,BeterminSandier2018,Lauritsen2021,marzo2025improvedlowerboundlogarithmic}, where the aforementioned energies are studied from both theoretical and numerical perspectives. Notably, Smale's 7th problem from his celebrated list \cite{Smale1998}, which is wide open and was first posed by Shub and Smale in \cite{ShubSmale1993}, concerns precisely the generation of point configurations on the usual $2$-sphere whose logarithmic energy is nearly minimal. In this setting, the actual minimizers are sometimes referred to as \emph{elliptic Fekete points}.

In recent years, interest in point configurations in spaces other than the sphere has grown significantly. Much of this attention has focused on projective spaces \cite{Pedro2023,Andersonetal2023,Lizarte,BFLG2025}, flat tori \cite{Marzo2018,Borda2024}, the special orthogonal group \cite{Mitchell,BeltranFerizovic2020}, and Grassmannians \cite{Grass2023,Pedro2025}. The problem can also be studied in a more general setting, where the space is any compact subset of $\R^n$. In this context, some celebrated results such as the poppy seed bagel theorem have been obtained; see \cite{BorodachovHardinSaff2019} for a comprehensive monograph.

In this paper, we focus on the special orthogonal group $\so\subseteq\R^{3\times 3}$, that is, the group of rotations in $\R^3$, consisting of all $3\times 3$ orthogonal matrices with determinant equal to $1$. Our goal is to find collections of matrices $O_1,\ldots,O_n\in\so$ such that their logarithmic energy is as low as possible. This logarithmic energy is given by
\begin{align*}
	\E(O_1,\ldots,O_n)&\defined\sum_{i\neq j}\log\norm{O_i-O_j}_F^{-1}\\
	&=-\frac{1}{2}\sum_{i\neq j}\log\norm{O_i-O_j}_F^2\\
	&=-\frac{1}{2}\sum_{i\neq j}\log(6-2\trace(O_i\transpose O_j)),
\end{align*}
where $\norm{O}_F=(\trace(O\transpose O))^{1/2}$ is the Frobenius norm of the matrix $O$. The main problems we address in this work are the following:
\begin{enumerate}
\item What is the lowest possible value of $\E(O_1,\ldots,O_n)$? From the general theory \cite{BorodachovHardinSaff2019} and the results in \cite{Andersonetal2023}, we know that
\begin{equation}\label{eq:minimo}
	\min_{O_1,\ldots,O_n\in\so}\E(O_1,\ldots,O_n)=\kappa n^2-\frac13n\log n+O(n),
\end{equation}
where
\[
\kappa=\frac{1}{\vol(\so)^2}\int_{O,O'\in\so}\log\norm{O-O'}_F^{-1}\dd O\dd O'
\]
is the continuous energy of the uniform measure, that is, the expected logarithmic energy of two matrices $O,O'\in\so$ chosen independently at random according to the standard Haar measure on $\so$. We will see below (see \cref{prop:continuousso3}) that
\[
\kappa=-\frac{1+\log2}{2}.
\]

\item Can we generate collections $O_1,\ldots,O_n$ whose energy is as close as possible to this minimal value? 
\end{enumerate}

Hence, we aim to find collections of special orthogonal matrices for which the $O(n)$ term in \eqref{eq:minimo} is provably as low as possible. To the best of our knowledge, there is currently only one known sequence of point sets on $\so$ for which some energy has been theoretically computed: the harmonic ensemble on $\so$, studied in \cite{BeltranFerizovic2020}. In that work, the logarithmic energy is not computed; instead, the authors study the Green energy and the Riesz $s$-energy for $s\in(0,3)$, defined respectively as
\begin{align*}
	\mathcal{E}_{\mathcal{G}}(O_1,\ldots,O_n)&\defined\sum_{i\neq j}\mathcal G(O_i,O_j),\\
	\mathcal{E}_{s}(O_1,\ldots,O_n)&\defined\sum_{i\neq j}\norm{O_i-O_j}_F^{-s},
\end{align*}
where $\mathcal G$ is the Green function of the Laplace--Beltrami operator on $\so$. For each $n$ of the form
\begin{equation}\label{eq:nL}
n=\binom{2L+3}{3}=\frac43L^3+O(L^2),\qquad L=1,2,\dotsc,
\end{equation}
the harmonic ensemble on $\SO(3)$ is a collection of $n$ random matrices $O_1,\ldots,O_n$ with a certain distribution on $\so\times\dotsb\times\so$ ($n$ times). The expected logarithmic energy of such configurations can be studied using the same techniques as in \cite{BeltranFerizovic2020}, as we do in \cref{sec:energy-harmonic}. Alternatively, it can be derived from the results in \cite{Andersonetal2023}, using the fact that $\so$ is essentially isometric to the real projective space $\RP^3$. 

\begin{theorem}\label{thm:energy-harmonic}
	The expected logarithmic energy of the points coming from the harmonic ensemble on $\SO(3)$ is, for $n=10,35,84,165,286,455,680,969,1330,1771,\ldots$ of the form \eqref{eq:nL},
	\begin{equation*}
		\kappa n^2-\frac{1}{3}n\log n+Cn+o(n),
	\end{equation*}
	with
	\begin{equation*}
		C=\frac{7}{3}-\gamma+\frac{\log 2}{6}-\frac{\log 3}{3}\approx 1.5054,
	\end{equation*}
	where $\gamma$ is the Euler--Mascheroni constant.
\end{theorem}

Comparing the energy in \cref{thm:energy-harmonic} with the theoretical minimum in \eqref{eq:minimo}, we observe that the harmonic ensemble on $\so$ matches the first two asymptotic terms. 

The main result of this paper, which represents a significant improvement over \cref{thm:energy-harmonic}, can be stated as follows.

\begin{theorem}[Main result]\label{th:main}
	For all $n$ in a certain infinite subsequence of the natural numbers \embparen{starting $n=2,3,8,10,12,21,24,27,30,33,36,52,56,\ldots$}, there exists a randomized constructive way to produce $n$ rotation matrices $O_1,\dotsc,O_n\in\SO(3)$ such that the expected value of the logarithmic energy of $O_1,\dotsc,O_n$ is
	\begin{equation*}
		\kappa n^2-\frac13n\log n + C_{\textup{zeros}}n+o(n),
	\end{equation*}
	where $C_{\textup{zeros}}=0.9190058\dotso$ is a constant.
\end{theorem}

See \cref{prop:zeros} for the actual construction and \eqref{eq:czeros} for a precise description of the constant $C_{\textup{zeros}}$.

The main tool we use is \cref{th:logenergy}, which establishes that, given any choice of $r$ spherical points, the logarithmic energy of the associated collection of $n=rs$ rotation matrices (obtained by choosing $s$ equidistant rotation matrices that fix each point) can be computed by studying another energy that depends solely on the spherical points. This allows us to study this other energy for different collections of spherical points. The configuration yielding \cref{th:main} is described in \cref{subsec:zeros}.

\subsection*{Notation}

The following notations will be used throughout this work. Let $f$ and $g$ be two real-valued functions, where we assume that $g(x)$ is nonzero for sufficiently large $x$.
\begin{itemize}
	\item The expression $f(x)\lesssim g(x)$ means that there exists a constant $C>0$ such that $f(x)\leq Cg(x)$ for all $x$.
	
	\item We say that $f(x)=O(g(x))$ as $x\tendsto\infty$ if there exists a positive constant $C$ and $x_0\in \R$ such that $\abs{f(x)}\leq C\abs{g(x)}$ for all $x>x_0$. We can express this condition as
	\begin{equation*}
		\limsup_{x\tendsto\infty} \abs[\bigg]{\frac{f(x)}{g(x)}}\leq C.
	\end{equation*}
	
	\item The expression $f(x)\asymp g(x)$ as $x\tendsto\infty$ means that $f(x)=O(g(x))$ and $g(x)=O(f(x))$, that is, $f(x)$ and $g(x)$ are asymptotically of the same order.
	
	\item We say that $f(x)=o(g(x))$ as $x\tendsto\infty$ if
	\begin{equation*}
		\lim_{x\tendsto\infty}\abs[\bigg]{\frac{f(x)}{g(x)}}=0.
	\end{equation*}
	
	\item Finally, $f(x)\sim g(x)$ as $x\tendsto\infty$  means that 
	\begin{equation*}
		\lim_{x\tendsto\infty}\frac{f(x)}{g(x)}=1.
	\end{equation*}
\end{itemize}
We have defined this asymptotic notation for $x\tendsto\infty$, but one can substitute $\infty$ with any real number $a$ and the definitions hold in all cases considering neighborhoods of $a$ when necessary.

\subsection*{Structure of the paper}

This work is structured as follows. In \cref{sec:cont-energy}, we compute $\kappa$, that is, the continuous logarithmic energy of the uniform measure on $\SO(3)$. Next, in \cref{sec:energy-harmonic}, we estimate the expected logarithmic energy of the points coming from the harmonic ensemble on $\SO(3)$, as described in \cite{BeltranFerizovic2020}. We then introduce, in \cref{sec:construction}, our general randomized construction for generating points on $\SO(3)$ from collections of points on $\S^2$. The expected energy of this construction is computed in \cref{sec:energy-our-construction}. In \cref{sec:colecciones}, we explore several specific instances of our construction, each derived from a collection of well-distributed points on the sphere, and compute their expected energies, including the case that leads to our main result in \cref{th:main}. In \cref{sec:auxiliary-results}, we include the proofs of some auxiliary results used throughout the paper. Finally, \cref{appendix:orthogonal polynomials} collects various definitions and results concerning orthogonal polynomials and special functions.

\section{The continuous logarithmic energy}\label{sec:cont-energy}

In this section, we compute the continuous energy $\kappa$, that is, the expected value of the energy for random, uniformly chosen  $O,O'\in\so$.

\begin{proposition}\label{prop:continuousso3}
	Let $O,O'\in\so$ be chosen uniformly and independently with respect to the Haar measure on $\so$. Then, we have
	\begin{align*}
		\kappa=\expectation[\big]{O,O'\in\so}{\log\norm{O-O'}_F^{-1}}=-\frac{1+\log2}2.
	\end{align*}
\end{proposition}

\begin{proof}
	Using a symmetry argument, we have
	\begin{equation*}
		\expectation[\big]{O,O'\in\so}{\log\norm{O-O'}_F^{-1}}=\expectation[\big]{O\in\so}{\log\norm{O-\Id_3}_F^{-1}},
	\end{equation*}
	where $\Id_3$ is the identity matrix of size $3$. This last expected value can be computed using \cite[Eq. (1) and Lemma 3.1]{BeltranFerizovic2020}, which yield
	\begin{align*}
		\expectation[\big]{O\in\so}{\log\norm{O-\Id_3}_F^{-1}}&=-\frac2\pi\int_0^\pi\log\biggl(\sqrt{8}\sin\frac{\theta}{2}\biggr)\sin^2\frac{\theta}2\dd\theta\\
		&=-\frac2\pi\int_0^{\pi/2}\log(8\sin^2\alpha)\sin^2\alpha\dd\alpha\\
		&=-\frac{\log8}{2}-\frac4\pi\int_0^{\pi/2}\log(\sin\alpha)\sin^2\alpha\dd\alpha\\
		&=-\frac{\log8}{2}-\frac12(1-\log4)\\
		&=-\frac{1+\log2}2,
	\end{align*}
	where we have used \cite[4.384-9]{integrales}.
\end{proof}

\section{Expected logarithmic energy of the harmonic ensemble on $\so$}\label{sec:energy-harmonic}

In this section, we give a proof of \cref{thm:energy-harmonic}. To prove that result, we need to perform a detailed asymptotic analysis of certain integrals involving orthogonal polynomials. For the reader’s convenience, we have encapsulated the necessary technical results as lemmas and deferred the proofs to \cref{sec:auxiliary-results}.

\begin{lemma}\label{lemma:aux-integral-1}
	The following equality holds:
	\begin{equation*}
		\int_{0}^{L^{-1/2}}\Gegenbauer{2L}{2}(\cos\theta)^2\sin^2\theta\dd\theta=\frac{\pi}{3}L^3+O(L^{5/2}).
	\end{equation*}
\end{lemma}

\begin{proof}
	See \cref{sec:proof-lemma:aux-integral-1}
\end{proof}

\begin{lemma}\label{lemma:aux-integral-2}
	The following equality holds:
	\begin{equation*}
		\int_{L^{-1/2}}^{\pi/2}\Gegenbauer{2L}{2}(\cos\theta)^2\sin^2\theta\dd\theta=O(L^{5/2}).
	\end{equation*}
\end{lemma}

\begin{proof}
	See \cref{sec:proof-lemma:aux-integral-2}
\end{proof}

\begin{lemma}\label{lemma:aux-integral-3}
	The following equality holds:
	\begin{equation*}
		\int_{L^{-1/2}}^{\pi/2}\log(\sin\theta)\Gegenbauer{2L}{2}(\cos\theta)^2\sin^2\theta\dd\theta=O(L^{5/2}\log L).
	\end{equation*}
\end{lemma}

\begin{proof}
	See \cref{sec:proof-lemma:aux-integral-3}
\end{proof}

\begin{lemma}\label{lemma:aux-integral-4}
	The following equality holds:
	\begin{align*}
		\MoveEqLeft\int_{0}^{L^{-1/2}}\log(\sin\theta)\Gegenbauer{2L}{2}(\cos\theta)^2\sin^2\theta\dd\theta\\
		&=-\frac{\pi}{3}L^3\log L+\biggl(\frac{(7-3\gamma-6\log{2})\pi}{9}\biggr)L^3+o(L^3).
	\end{align*}
\end{lemma}

\begin{proof}
	See \cref{sec:proof-lemma:aux-integral-4}
\end{proof}

\begin{proof}[Proof of \cref{thm:energy-harmonic}]
	Let $\mathbb{E}=\expectation{x\in(\SO(3))^n}{\E(x)}$ be the expected logarithmic energy. Then, following the first lines of the proof of \cite[Proposition 3.2]{BeltranFerizovic2020},
	\begin{align*}    
		\mathbb{E}&=\kappa n^2 +\frac{2}{\pi}\int_{0}^{\pi}\log\biggl(\sqrt{8}\sin\frac{t}{2}\biggr)\Gegenbauer{2L}{2}\biggl(\cos\frac{t}{2}\biggr)^2\sin^2\frac{t}{2}\dd t\\
		&=\kappa n^2 +\frac{4}{\pi}\int_{0}^{\pi/2}\log\bigl(\sqrt{8}\sin\theta\bigr)\Gegenbauer{2L}{2}(\cos\theta)^2\sin^2\theta\dd\theta.
	\end{align*}
	Let
	\begin{equation*}
		I=\frac{4}{\pi}\int_{0}^{\pi/2}\log\bigl(\sqrt{8}\sin\theta\bigr)\Gegenbauer{2L}{2}(\cos\theta)^2\sin^2\theta\dd\theta.
	\end{equation*}
	We have $I=I_1+I_2$, where
	\begin{align*}
		I_1&=\frac{6\log 2}{\pi}\int_{0}^{\pi/2}\Gegenbauer{2L}{2}(\cos\theta)^2\sin^2\theta\dd\theta,\\
		I_2&=\frac{4}{\pi}\int_{0}^{\pi/2}\log(\sin\theta)\Gegenbauer{2L}{2}(\cos\theta)^2\sin^2\theta\dd\theta.
	\end{align*}
	Using \cref{lemma:aux-integral-1,lemma:aux-integral-2}, we have
	\begin{align*}
		I_1
		=&\frac{6\log 2}{\pi}\biggl(\frac{\pi}{3}L^3+O(L^{5/2})\biggr)=2L^3\log2+o(L^3).
	\end{align*}
	Using \cref{lemma:aux-integral-3,lemma:aux-integral-4}, we have
	\begin{align*}
		I_2=&\frac{4}{\pi}\biggl(-\frac{\pi}{3}L^3\log L+\biggl(\frac{(7-3\gamma-6\log{2})\pi}{9}\biggr)L^3+o(L^{3})\biggr)\\
		=&-\frac{4}{3}L^3\log L+\biggl(\frac{7-3\gamma-6\log 2}{3}\biggr)\frac{4}{3}L^3+o(L^{3}).
	\end{align*}
	Then, the expected energy is
	\begin{align*}    
		\mathbb{E}=\kappa n^2 -\frac{4}{3}L^3\log L+\biggl(\frac{7-3\gamma-6\log 2}{3}+\frac{3}{2}\log 2\biggr)\frac{4}{3}L^3+o(L^{3}).
	\end{align*}
	The result follows using that $n=(4/3)L^3+O(L^2)$.
\end{proof}

\section{A construction on $\so$ from a set of points on $\s$}\label{sec:construction}

In this section, we describe how to obtain a constructive collection of points on $\so$ from any collection of points on $\s$.  A similar idea has been studied in \cite{Mitchell}, mainly from a numerical point of view. Let 
\begin{equation}\label{eq:Rphi}
	R(\phi)\defined\begin{pmatrix}
		\cos\phi& -\sin\phi & 0\\
		\sin\phi& \cos\phi & 0\\
		0 &0 & 1
	\end{pmatrix}
\end{equation}
be the rotation matrix with angle $\phi$ and with base point $e_3=(0,0,1)\transpose $.

Let $p\in\s$ be a spherical point. Let $s$ be a positive integer. We say that $s$ matrices $O_1,\ldots,O_s\in\so$ are \emph{equally spaced with base point $p$} if there exist an orthogonal matrix $H_p\in \so$ with $H_pe_3=p$ and an angle $\phi\in[0,2\pi)$ such that
\begin{equation*}
	O_j=H_p R(\phi_j+\phi).
\end{equation*}
Here, $\phi_j=2\pi j/s$, with $j=1,\dotsc,s$, are the $s$ roots of unity. For example, the matrices $R(\phi_j)$, with $j=1,\ldots,s$, are equally spaced with base point $e_3$. Moreover, they are, for any fixed $s$, the unique matrices with that property up to the choice of a random phase.

One can give explicit formulas for $H_p$ in different manners. Here is a particular choice:
\begin{itemize}
\item If $p=e_3$, let $H_p^{(0)}=\diag(1,1,1)=\Id_3$, that is, the identity matrix.
\item If $p=-e_3$, let $H_p^{(0)}=\diag(1,-1,-1)$.
\item Otherwise, let $p=(x,y,z)\transpose $, with $x^2+y^2\neq0$, and let
\[
H_p^{(0)}\defined\begin{pmatrix}
	\frac{y}{\sqrt{x^2+y^2}} &\frac{zx}{\sqrt{x^2+y^2}}  & x\\
	\frac{-x}{\sqrt{x^2+y^2}} &\frac{zy}{\sqrt{x^2+y^2}} & y\\
0 &	-\sqrt{x^2+y^2} & z
\end{pmatrix}\eqcolon(p^2,p^1,p)\in\so.
\]
\end{itemize}

\begin{definition}\label{def:pointsinso3}
Given $r$ points $p_1,\dotsc,p_r\in\s$ and a positive integer $s$, the \emph{collection of special orthogonal matrices associated with $p_1,\ldots,p_r$ and $s$} is a random set of $rs$ elements on $\so$ obtained as follows:
\begin{enumerate}
	\item Choose $\phi^{(1)},\ldots,\phi^{(r)}$ uniformly at random in $[0,2\pi)$.
	
	\item Consider the set
	\[
	\bigcup_{i=1}^r\bigcup_{j=1}^s\set{O_{i,j}},\quad\text{where}\quad O_{i,j}=H_{p_i}^{(0)}R(\phi_j+\phi^{(i)}).
	\]
	That is, for each $p_i$, with $1\leq i\leq r$, we choose $s$ equally spaced matrices with base point $p_i$ rotated by a random phase $\phi^{(i)}$.
\end{enumerate}
Since the set described above is actually a random set of points, that is, a random process on $\so$, in order to study its properties we will compute the expected value of the logarithmic energy of the point set when the phases are chosen uniformly at random in $[0,2\pi)$.
\end{definition}

\section{The expected logarithmic energy of the points in the construction}\label{sec:energy-our-construction}
 
In this section, we prove the following result, which relates the logarithmic energy of the collection of $n$ special orthogonal matrices associated with a collection of spherical points to a certain expresion involving only the dot products of these spherical points (alternatively, we might write it in terms of the mutual distances between the spherical points). This connection enables us to compute the logarithmic energy of these sets of orthogonal matrices using only computations on the sphere, a topic that has been much more studied. 

\begin{theorem}\label{th:logenergy}
The expected logarithmic energy of the collection of special orthogonal matrices associated with $p_1,\ldots,p_r$ and $s$ is equal to
\begin{align*}
\MoveEqLeft\expectation[\bigg]{}{-\sum_{(i,j)\neq(i',j')}\log\|O_{i,j}-O_{i',j'}\|_F}\\
&=-\frac{n^2}{2}\log2+\frac{n}2\log2-n\log s-s^2\sum_{i\neq i'}\log\biggl(1+\sqrt{\frac{1-\innerprod{p_i}{p_{i'}}}{2}}\biggr).
\end{align*}
\end{theorem}

To prove \cref{th:logenergy}, we will use the following lemmas.

\begin{lemma}\label{lem:equallyspaced}
Let $p\in\s$ and let $O_1,\dotsc,O_s\in\so$ be $s$ equally spaced matrices with base point $p$. Then,
\[
\sum_{i\neq j}\log\norm{O_i-O_j}_F=\frac{s(s-1)}{2}\log2+s\log s.
\]
\end{lemma}

\begin{proof}
Let $S$ be the sum in the lemma. Note that
\begin{align*}
	S&=\log\prod_{i\neq j}\norm[\big]{H_p R(\phi_i+\phi)-H_p R(\phi_j+\phi)}_F\\
	 &=\log\prod_{i\neq j}\norm[\big]{\Id_3-R(\phi_j+\phi)R(\phi_i+\phi)^{-1}}_F\\
	 &=\frac12\log\prod_{i\neq j}\norm[\big]{\Id_3-R(\phi_j-\phi_i)}_F^2\\
	 &=\frac12\log\prod_{i\neq j}\bigl(4-4\cos(\phi_j-\phi_i)\bigr)\\
	 &=s(s-1)\log2+\frac12\log\prod_{i=1}^s\prod_{\substack{j=1\\j\neq i}}^s\biggl(1-\cos\biggl(2\pi\frac{i-j}{s}\biggr)\biggr).
\end{align*}
The inner product can be rewritten, setting $k=i-j$, as
\begin{equation*}
\prod_{\substack{j=1\\j\neq i}}^s\biggl(1-\cos\biggl(2\pi\frac{i-j}{s}\biggr)\biggr)=
\prod_{k=1}^{i-1}\biggl(1-\cos\biggl(2\pi\frac{k}{s}\biggr)\biggr)
\prod_{k=i-s}^{-1}\biggl(1-\cos\biggl(2\pi\frac{k}{s}\biggr)\biggr).
\end{equation*}
This last part, taking $\hat k=k+s$ and using the periodicity of the cosine, equals
\begin{equation*}
	\prod_{k=i-s}^{-1}\biggl(1-\cos\biggl(2\pi\frac{k}{s}\biggr)\biggr)=
	\prod_{\hat k=i}^{s-1}\biggl(1-\cos\biggl(2\pi\frac{\hat k-s}{s}\biggr)\biggr)=\prod_{\hat k=i}^{s-1}\biggl(1-\cos\biggl(2\pi\frac{\hat k}{s}\biggr)\biggr).
\end{equation*}
Putting everything together, we have proved that
\begin{equation*}
	\prod_{\substack{j=1\\j\neq i}}^s\biggl(1-\cos\biggl(2\pi\frac{i-j}{s}\biggr)\biggr)=
	\prod_{k=1}^{s-1}\biggl(1-\cos\biggl(2\pi\frac{k}{s}\biggr)\biggr).
\end{equation*}
In particular, its value is independent of $i$. Hence, we have
\begin{align*}
	S&=s(s-1)\log2+\frac{s}2\log\prod_{j=1}^{s-1}\biggl(1-\cos\biggl(\frac{2\pi j}{s}\biggr)\biggr)\\
	&=s(s-1)\log2+\frac{s(s-1)}{2}\log2+s\log\prod_{j=1}^{s-1}\sin\biggl(\frac{\pi j}{s}\biggr),
\end{align*}
where we have used that $1-\cos(2x)=2\sin^2 x$. From \cite[1.392-1]{integrales}, we have
\[
\prod_{j=0}^{s-1}\sin\Big(x+\frac{\pi j}{s}\Big)=\frac{\sin(sx)}{2^{s-1}},
\]
and hence as $x\tendsto 0$
\[
\prod_{j=1}^{s-1}\sin\Big(x+\frac{j\pi}{s}\Big)=\frac{1}{2^{s-1}}\frac{\sin(sx)}{\sin(x)}\Longrightarrow \prod_{j=1}^{s-1}\sin\Big(\frac{j\pi}{s}\Big)=\frac{s}{2^{s-1}}.
\]
We thus have
\begin{equation*}
	S=\frac{3}{2}s(s-1)\log2+s\log\frac{s}{2^{s-1}},
\end{equation*}
and the lemma follows.
\end{proof}

\begin{lemma}\label{lem:an_integral}
	Let $H\in\so$. Then, for $\alpha,\beta\in\R$ with $|\beta|\leq\alpha/3$, we have
	\[
	\frac1{2\pi}\int_0^{2\pi}\log\Bigl(\alpha+\beta\trace\bigl(HR(\phi)\bigr)\Bigr)\dd\phi= 2\log\frac{\sqrt{\alpha-\beta}+\sqrt{\alpha+\beta+2\beta e_3\transpose He_3}}{2}.
	\]
	In particular, taking $\alpha=6$ and $\beta=-2$, we have
	\[
	\frac1{2\pi}\int_0^{2\pi}\log\Bigl(6-2\trace\bigl(HR(\phi)\bigr)\Bigr)\dd\phi=2\log\Bigl(\sqrt{2}+\sqrt{1-e_3\transpose He_3}\Bigr).
	\]
\end{lemma}
\begin{proof}
Let $I$ be the integral in the lemma. Writing $h_{ij}$ for the entries of $H$, we have
\begin{align*}
	I&=\frac1{2\pi}\int_0^{2\pi}\log\Bigl(\alpha+\beta\bigl((h_{11}+h_{22})\cos\phi+(h_{12}-h_{21})\sin\phi+h_{33}\bigr)\Bigr)\dd\phi\\
	&=\log A+\frac1{2\pi}\int_0^{2\pi}\log\frac{\alpha+\beta\bigl((h_{11}+h_{22})\cos\phi+(h_{12}-h_{21})\sin\phi+h_{33}\bigr)}{A}\dd\phi,
\end{align*}
for any choice of $A>0$. Choosing $A=\alpha+\beta h_{33}$ and defining $a=\beta(h_{12}-h_{21})/A$ and $b=\beta(h_{11}+h_{22})/A$, we get
\begin{equation*}
	I=\log (\alpha+\beta h_{33})+\frac1{2\pi}\int_0^{2\pi}\log(1+a\sin\phi+b\cos\phi)\dd\phi.
\end{equation*}
Note that
\begin{equation*}
	a^2+b^2=\frac{h_{11}^2+h_{22}^2+h_{21}^2+h_{12}^2+2(h_{11}h_{22}-h_{12}h_{21})}{(\alpha+\beta h_{33})^2}\beta^2.
\end{equation*}
From the fact that the columns of $H$ have norm $1$, we get
\begin{equation*}
	h_{11}^2+h_{22}^2+h_{21}^2+h_{12}^2=2-h_{31}^2-h_{32}^2,
\end{equation*}
which in turn, using that the last row of $H$ also has norm $1$, equals $1+h_{33}^2$. Moreover, since the determinant of $H$ is equal to $1$, from the construction of the inverse we have that the $(3,3)$ entry of $H^{-1}$ is precisely $h_{11}h_{22}-h_{12}h_{21}$, and from the fact that $H^{-1}=H\transpose $ we deduce that this last expresion equals $h_{33}$. Putting everything together, we have
\begin{equation*}
	a^2+b^2=\frac{1+h_{33}^2+2h_{33}}{(\alpha+\beta h_{33})^2}\beta^2=\frac{(1+h_{33})^2}{(\alpha+\beta h_{33})^2}\beta^2\leq1,
\end{equation*}
where the last inequality follows from the facts that $\abs{h_ {33}}\leq 1$ and $\abs{\beta}\leq\alpha/3$. Hence, we may use \cite[(4.225,3)]{integrales}, which yields
\begin{align*}
	I&=\log (\alpha+\beta h_{33})+\log\frac{1+\sqrt{1-\beta^2(1+h_{33}^2+2h_{33})/(\alpha+\beta h_{33})^2}}{2}\\
	&=-\log2+\log\Bigl(\alpha+\beta h_{33}+\sqrt{(\alpha+\beta h_{33})^2-(1+h_{33}^2+2h_{33})\beta^2}\Bigr)\\
	&=-\log2+\log\Bigl(\alpha+\beta h_{33}+\sqrt{\alpha- \beta}\sqrt{\alpha+\beta+2\beta h_{33}}\Bigr)\\
	&=-\log2+\log\biggl(\frac12\bigl(\sqrt{\alpha-\beta}+\sqrt{\alpha+\beta+2\beta h_{33}}\bigr)^2\biggr)\\
	&=-\log4+2\log\Bigl(\sqrt{\alpha-\beta}+\sqrt{\alpha+\beta+2\beta h_{33}}\Bigr).
\end{align*}
Since $h_{33}=e_3\transpose He_3$, the lemma follows.
\end{proof}

\begin{lemma}\label{lem:twopoints}
	Let $p,p'\in\s$ and let $\set{O_j}_{j=1}^{s}$ and $\set{O'_j}_{j=1}^{s}$ be two collections of $s$ equally spaced matrices with respective base points $p$ and $p'$, that is,
	\begin{align*}
		O_j&=H_pR(\phi_j+\phi),\\
		O'_j&=H_{p'}R(\phi_j+\phi'),
	\end{align*}
	where $\phi,\phi'\in[0,2\pi)$ are chosen uniformly at random, and $H_p,H_{p'}\in\SO(3)$ are orthogonal matrices such that $H_pe_3=p$ and $H_{p'}e_3=p'$. Then, the expected value of the crossed logarithmic energy is
	\[
	\expectation[\bigg]{}{\sum_{i=1}^s\sum_{j=1}^s\log\norm{O_i-O'_j}_F}=s^2\log\bigl(\sqrt2+\sqrt{1-\innerprod{p}{p'}}\bigr).
	\]
 \end{lemma}
 
\begin{proof}
	Let $S$ be the expected value in the lemma. We have
	\begin{align*}
		S&=\frac{1/2}{(2\pi)^2}\int_0^{2\pi}\int_0^{2\pi}\sum_{i=1}^s\sum_{j=1}^s\log\norm[\big]{H_pR(\phi_i+\phi)-H_{p'}R(\phi_j+\phi')}_F^2\dd\phi\dd\phi'\\
		&=\sum_{i,j=1}^s\frac{1/2}{(2\pi)^2}\int_0^{2\pi}\int_0^{2\pi}\log\Bigl(6-2\trace\bigl(R(-\phi_j-\phi')H_{p'}\transpose H_pR(\phi_i+\phi)\bigr)\Bigr)\dd\phi\dd\phi'\\
		&=\sum_{i,j=1}^s\frac{1/2}{(2\pi)^2}\int_0^{2\pi}\int_0^{2\pi}\log\Bigl(6-2\trace\bigl(H_{p'}\transpose H_pR(\phi_i-\phi_j+\phi-\phi')\bigr)\Bigr)\dd\phi\dd\phi',
	\end{align*}
	where in the last equality we have used the cyclic property of the trace. Due to the periodicity of $R(\phi)$ as defined in \eqref{eq:Rphi}, for any fixed $\phi_i,\phi_j$, and $\phi'$ we can change the integration interval of the variable $\phi$ from $[0,2\pi]$ to $[\phi_i+\phi_j+\phi',\phi_i+\phi_j+\phi'+2\pi]$ without changing the value of the integral. Thus, we get
	\begin{align*}
		S&=\sum_{i,j=1}^s\frac{1/2}{(2\pi)^2}\int_0^{2\pi}\int_0^{2\pi}\log\Bigl(6-2\trace\bigl(H_{p'}\transpose H_pR(\phi)\bigr)\Bigr)\dd\phi\dd\phi'\\
		&=\frac{s^2}{4\pi}\int_0^{2\pi}\log\Bigl(6-2\trace\bigl(H_{p'}\transpose H_pR(\phi)\bigr)\Bigr)\dd\phi.
	\end{align*}
	Finally, from \cref{lem:an_integral} we conclude that
	\begin{equation*}
		S=s^2\log\Bigl(\sqrt2+\sqrt{1-e_3\transpose H_{p'}\transpose H_pe_3}\Bigr).\\
	\end{equation*}
	Since $e_3\transpose H_{p'}\transpose H_pe_3=p'{\transpose} p=\innerprod{p}{p'}$, the lemma follows.
\end{proof}

We are finally ready to prove \cref{th:logenergy}.

\begin{proof}[Proof of \cref{th:logenergy}]
	Let $S=S_1+S_2$ be the expected value in the lemma, where
	\begin{align*}
		S_1&=\expectation[\bigg]{}{-\sum_{i=1}^r\sum_{j\neq j'}\log\norm{O_{i,j}-O_{i,j'}}_F},\\
		S_2&=\expectation[\bigg]{}{-\sum_{i\neq i'}\sum_{j,j'=1}^s\log\norm{O_{i,j}-O_{i',j'}}_F}.
	\end{align*}
	From Lemma \ref{lem:equallyspaced}, we have
	\[
	S_1=-r\biggl(\frac{s(s-1)}{2}\log2+s\log s\biggr).
	\]
	Moreover, from Lemma \ref{lem:twopoints}, we have
	\[
	S_2=-\sum_{i\neq i'}s^2\log\bigl(\sqrt2+\sqrt{1-\innerprod{p_i}{p_{i'}}}\bigr).
	\]
	Adding the two terms and simplifying finishes the proof.
\end{proof}

\section{Collections of points on $\s$ and the associated energy on $\so$}\label{sec:colecciones}

With Theorem \ref{th:logenergy} in hand, we can now choose collections of points on $\s$ and compute the logarithmic energy of the associated collections of special orthogonal matrices. The expresion
\begin{equation}\label{eq:energiaesfera}
	\sum_{i\neq i'}\log\biggl(1+\sqrt{\frac{1-\innerprod{p_i}{p_{i'}}}{2}}\biggr)
\end{equation}
in Theorem \ref{th:logenergy} can also be seen as some kind of discrete pairwise interaction energy on $\s$. Unlike other classical potentials, however, this energy remains finite for all configurations and does not diverge as points approach each other. Let us start by computing the continuous energy, that is, the expected value of this energy for random, uniformly chosen $p,p'\in\s$.

\begin{proposition}\label{prop:continuousesfera}
Let $p,p'\in\s$ be chosen uniformly and independently. Then, we have
\begin{equation*}
	\expectation[\Bigg]{p,p'\in\s}{\log\biggl(1+\sqrt{\frac{1-\innerprod{p}{p'}}2}\biggr)}=\frac12.
\end{equation*}
\end{proposition}

\begin{proof}
	By rotational symmetry, we may fix $p'=e_3=(0,0,1)\transpose$. Then, we have
	\begin{equation*}
		\expectation[\Bigg]{p,p'\in\s}{\log\biggl(1+\sqrt{\frac{1-\innerprod{p}{p'}}2}\biggr)}=
		\expectation[\Bigg]{p\in\s}{\log\biggl(1+\sqrt{\frac{1-\innerprod{p}{e_3}}2}\biggr)}.
	\end{equation*}
	Since now the expected value depends only on $p_3$, we can express it as a univariate integral (see, for example, \cite[Lemma 1]{BeltranMarzoOrtega}) that can be easily computed using the change of variables $t=1-2u^2$:
	\begin{align*}
		\expectation[\Bigg]{p\in\s}{\log\biggl(1+\sqrt{\frac{1-\innerprod{p}{e_3}}2}\biggr)}&=\frac12\int_{-1}^1\log\biggl(1+\sqrt{\frac{1-t}{2}}\biggr)\dd t\\
		&=\frac12\int_0^14u\log(1+u)\dd u\\
		&=\frac12,
	\end{align*}
	as claimed. 
\end{proof}

A natural strategy to search for points on $\s$ that maximize \eqref{eq:energiaesfera}---and hence, from Theorem \ref{th:logenergy}, produce minimal values of the logarithmic energy on $\so$---is to consider configurations that are known to yield quasioptimal values for general energies, such as the logarithmic or Riesz energies, on the sphere. Such collections of points have been extensively studied, but theoretical bounds for the corresponding energies are known in just a few cases:

\begin{enumerate}
\item Random, independent, and identically distributed points on $\s$. \label{item:iid}

\item Zeros of random polynomials \cite{ArmentanoBeltranShub2011,delatorre2022expectedenergyzeroselliptic}: A random polynomial of degree $r$ is generated by drawing the coefficient of $z^j$ from a complex Gaussian distribution with mean $0$ and variance $\binom{r}{j}$. Then, its $r$ complex zeros are mapped to the sphere via the inverse stereographic projection.\label{item:zeros}

\item Equal area partitions \cite{Leopardi}: The sphere is partitioned into $r$ regions of equal area, each with diameter of the order of $O(1/\sqrt{r})$. A point is then randomly chosen within each region.\label{item:EAP}

\item The spherical ensemble \cite{AlishahiZamani}: Two random $r\times r$ matrices $A$ and $B$ are generated by drawing their entries from a complex Gaussian distribution with mean $0$ and variance $1$. Then, the $r$ complex eigenvalues of $A^{-1}B$ are computed and sent to $\s$ using the inverse stereographic projection.\label{item:sph-ensemble}

\item The harmonic ensemble \cite{BeltranMarzoOrtega}: A determinantal point process given by the spherical harmonics of bounded degree. Expected values of energies, discrepancy, and other quantities can be computed, although in practice the actual generation of points is more involved.\label{item:harmonic}

\item The diamond ensemble \cite{BeltranEtayo2020}: A number of parallels is chosen on the sphere. Then, on each parallel we place a number of equally spaced points rotated by a random phase. The heights of the parallels and the number of points per parallel must be carefully chosen.\label{item:diamond}

\item Zeros of random instances of the polynomial eigenvalue problem \cite{armentano2024logarithmicenergysolutionspolynomial}: A random instance of a polynomial eigenvalue problem with a certain distribution of the coefficients is generated. The zeros are again sent to the sphere via the inverse stereographic projection. This process generalizes both the construction based on the zeros of random polynomials and the one arising from the spherical ensemble.\label{item:EVP}
\end{enumerate}

Note that all these collections of points are randomized, in the sense that each construction depends on the choice of a number of random parameters. To date, such randomness has been essential if we demand that the (expected values of the) energies can be analytically computed. Among these collections, the diamond ensemble yields the best known expected values for the logarithmic energy; however, in this case, only the logarithmic energy (and not the Riesz energies) is known. The lowest provable expected Riesz energy is achieved by the zeros of random polynomials. 

In this work we focus on cases \ref{item:iid} to \ref{item:harmonic}, thus including the case of zeros of random polynomials. The remaining cases---the diamond ensemble and the polynomial eigenvalue problem---are much more involved. We leave them as an open problem.

\subsection{Uniformity of minimizers}

Before analyzing each specific choice of spherical points, we will prove that the kernel
\begin{equation*}
	K(x,y)=-\log\biggl(1+\sqrt{\frac{1-\innerprod{x}{y}}2}\biggr),\qquad x,y\in\S^2,
\end{equation*}
is conditionally strictly positive definite. Then, from \cite[Theorem 6.2.1]{BorodachovHardinSaff2019}, the uniform measure on $\S^2$ is the unique minimizer for the associated $K$-energy. This justifies the convenience of the different sequences of spherical points described at the beginning of this section.

\begin{definition}
	A function $f\from I\to\R$ is \emph{absolutely monotone} on an interval $I$ if, for any $n\geq 0$ and $x\in I$, its $n$th derivative $f^{(n)}(x)$ exists\footnote{If $x$ is an endpoint of $I$, then $f^{(n)}(x)$ means the appropriate one-sided derivative.} and is nonnegative. The function $f$ is \emph{strictly absolutely monotone} if, in addition, $f^{(n)}(x)>0$ for all $x$ in the interior of $I$. 
	
	A function $f\from I\to\R$ is \emph{completely monotone} on an interval $I$ if, for any $n\geq 0$ and $x\in I$, its $n$th derivative $f^{(n)}$ exists and satisfies $(-1)^nf^{(n)}(x)\geq0$. The function $f$ is \emph{strictly completely monotone} if, in addition, $(-1)^nf^{(n)}(x)>0$ for all $x$ in the interior of $I$.
\end{definition}

\begin{theorem}[{\cite[Theorem 4.4.7]{BorodachovHardinSaff2019}}]\label{th:447book}
If \( g\from (0, \infty) \to \R \) is a differentiable function such that \( -g'(t) \) is strictly completely monotone on \( (0, \infty) \) with \( \lim_{t \tendsto \infty} g'(t) = 0 \) and  
\[
g(0)\defined \lim_{t \tendsto 0^+} g(t),
\]  
then the kernel  
\[
K(x, y) = g\bigl( \norm{ x - y }^2 \bigr)
\]  
is conditionally strictly positive definite on \( A \times A \)  
for any compact set \( A \contained \mathbb{R}^p \)  
with nonzero \( K \)-capacity.

\end{theorem}

\begin{remark}\label{remark:trivial}
	It is trivial to see that $f(t)$ being strictly absolutely monotone is equivalent to $g(t)=f(1-t/2)$ being strictly completely monotone, and that, if $g(t)$ is strictly completely monotone, then so is $-g'(t)$.
\end{remark}

\begin{theorem}
	Let $x,y\in\S^2$. The kernel
	\begin{equation*}
		K(x,y)=-\log\biggl(1+\sqrt{\frac{1-\innerprod{x}{y}}2}\biggr)=-\log\biggl(1+\sqrt{\frac{\|x-y\|^2}4}\biggr)
	\end{equation*}
	is conditionally strictly positive definite.
\end{theorem}

\begin{proof}
	From \cref{th:447book} and \cref{remark:trivial}, with $g(t)=-\log\bigl(1+\sqrt{t/4}\bigr)$ and $f(t)=-\log\bigl(1+\sqrt{(1-t)/2}\bigr)$, it suffices to check that
	\begin{equation*}
		f(t)=-\log\biggl(1+\sqrt{\frac{1-t}{2}}\biggr)
	\end{equation*}
	is strictly absolutely monotone on $[-1,1]$. To compute the $n$th derivative of $f$ we use the Faà di Bruno's formula, which generalizes the chain rule to higher derivatives. Consider the functions $g\from[0,\sqrt{2}]\to\R$ and $h\from[-1,1]\to\R$ given by
	\begin{align*}
		g(x)&=-\log(1+x/\sqrt{2}),\\
		h(x)&=\sqrt{1-x}.
	\end{align*}
	Note that $f=g\composed h$. Using the Faà di Bruno's formula, we have
	\begin{equation}\label{eq:nth-derivative-kernel}
		f^{(n)}(t)=\frac{\dd^n}{\dd t^n}g(h(t))
		=\sum_{k=1}^{n}g^{(k)}(h(t))
		B_{n,k}(h'(t),\dotsc,h^{(n-k+1)}(t)),
	\end{equation}
	where 
	\begin{equation*}
		B_{n,k}(h'(t),\dotsc,h^{(n-k+1)}(t))=\hspace{1cm}\sum_{\mathclap{\substack{j_1+\dotsb+j_{n-k+1}=k\\ j_1+2j_1+\dotsb+(n-k+1)j_{n-k+1}=n}}}
		\hspace{1cm}\frac{n!}{j_1!j_2!\dotsb j_{n-k+1}!}\prod_{\ell=1}^{n-k+1}\biggl(\frac{h^{(\ell)}(x)}{\ell!}\biggr)^{j_{\ell}}
	\end{equation*}
	is a Bell polynomial. Note that
	\begin{align*}
		g^{(k)}(x)&=(-1)^{k}(k-1)!(\sqrt{2}+x)^{-k},\\
		h^{(k)}(x)&=-\frac{(2k-3)!!}{2^k}(1-x)^{-(2k-1)/2}.
	\end{align*}
	To check that all the summands in \eqref{eq:nth-derivative-kernel} are positive, note that, for each $k$, with $1\leq k\leq n$, we have
	\begin{align*}
		\prod_{\ell=1}^{n-k+1}\biggl(\frac{h^{(\ell)}(x)}{\ell!}\biggr)^{j_{\ell}}&=\prod_{\ell=1}^{n-k+1}(-1)^{j_{\ell}}\Bigl(\frac{(2\ell-3)!!}{2^{\ell}\ell!}(1-x)^{-(2\ell-1)/2}\Bigr)^{j_{\ell}}\\
		&=(-1)^{j_{1}+\dotsb+j_{n-k+1}}\prod_{\ell=1}^{n-k+1}\Bigl(\frac{(2\ell-3)!!}{2^{\ell}\ell!}(1-x)^{-(2\ell-1)/2}\Bigr)^{j_{\ell}}\\
		&=(-1)^{k}\prod_{\ell=1}^{n-k+1}\Bigl(\frac{(2\ell-3)!!}{2^{\ell}\ell!}(1-x)^{-(2\ell-1)/2}\Bigr)^{j_{\ell}}.
	\end{align*}
	Therefore, we have
	\begin{align*}
		f^{(n)}(t)&=n!\sum_{k=1}^{n}(k-1)!(\sqrt{2}+\sqrt{1-t})^{-k}\\
		&\times\sum_{\substack{j_1+\dotsb+j_{n-k+1}=k\\ j_1+2j_1+\dotsb+(n-k+1)j_{n-k+1}=n}} \prod_{\ell=1}^{n-k+1}\frac{1}{j_{\ell}!}\Bigl(\frac{(2\ell-3)!!}{2^{\ell}\ell!}(1-t)^{-(2\ell-1)/2}\Bigr)^{j_{\ell}}.
	\end{align*}
	Since all the summands are positive, we conclude that $f^{(n)}(t)\geq 0$ for all $t\in[-1,1]$ and hence $f$ is strictly absolutely monotone on $[-1,1]$.
\end{proof}

\subsection{Point set \ref{item:iid}: random uniform points on $\s$}

We start by assuming that $p_1,\ldots,p_r$ are chosen independently and identically distributed (i.i.d.) with respect to the uniform distribution on $\s$.

\begin{proposition}\label{prop:iid}
Assume that $p_1,\ldots,p_r\in\s$ are i.i.d. on $\s$. Then, the expected logarithmic energy of the corresponding collection of $n=rs$ special orthogonal matrices associated with $p_1,\ldots,p_r$ and a positive integer $s$ is
\[
\mathbb{E}_{\unif}=\kappa n^2-n\log s+\frac{n}{2}s+ \frac{n}2\log2.
\]
\end{proposition}

It is immediate to see that, for a fixed value of $n=rs$, the minimum of this quantity is attained when $s=2$. This leads to our first complete construction of $n=2r$ poinst on $\so$: \emph{the collection of special orthogonal matrices associated with i.i.d. points $p_1,\ldots,p_r$ and $s=2$.} By comparison with \eqref{eq:minimo}, we observe that the logarithmic energy satisfies
\[
\mathbb{E}_{\unif}(\text{any $r$; $s=2$})=\kappa n^2+O(n),
\]
which matches the theoretical minimum up to the leading order term.

\begin{proof}[Proof of \cref{prop:iid}]

From \cref{th:logenergy} and \cref{prop:continuousesfera}, we deduce that
\begin{align*}
	\mathbb{E}_{\unif}&=-\frac{n^2}{2}\log2+\frac{n}2\log2-n\log s-\frac12s^2r(r-1)\\
	&=-\frac{n^2}{2}\log2+\frac{n}2\log2-n\log s-\frac12n(n-s).
\end{align*}
The proposition follows.
\end{proof}

\subsection{Point set \ref{item:zeros}: zeros of random polynomials}\label{subsec:zeros}
We now consider the case where the points $p_1,\ldots,p_r\in\S^2$ are chosen as the inverse stereographic projections of the zeros of elliptic random polynomials (sometimes called Kostlan--Shub--Smale polynomials or Bombieri polynomials), as described in \cite{ArmentanoBeltranShub2011}. That is, a random polynomial of degree $r$ is generated by drawing the coefficient of $z^j$ from a complex Gaussian distribution with mean $0$ and variance $\binom{r}{j}$. Then, its $r$ complex zeros are mapped to the sphere via the inverse stereographic projection.

\begin{proposition}\label{prop:zeros}
	Let $p_1,\dotsc,p_r\in\s$ be the inverse stereographic projections of the zeros of elliptic random polynomials. Then, the expected logarithmic energy of the corresponding collection of $n=rs$ special orthogonal matrices associated with $p_1,\ldots,p_r$ and a positive integer $s$ is given by
	\[
	\mathbb{E}_{\textup{zeros}}=\kappa n^2-n\log s+\frac{n}2\log2-\sqrt{n}s^{3/2}J+o\bigl(\sqrt{n}s^{3/2}\bigr),
	\]
	where $J=-0.5787893\ldots$ is as in \eqref{eq:J}.
\end{proposition}

It is straightforward to see that, for fixed $n=rs$, the minimum of this quantity is attained when $s=(4n/(9J^2))^{1/3}$, or, equivalently, $s=\sqrt{4r/(9J^2)}$. Since this last number is not an integer, we take its integer part to define our second complete construction of $n=r\floor{\sqrt{4r/(9J^2)}}$ points on $\so$: \emph{the collection of special orthogonal matrices associated with the inverse stereographic projections $p_1,\dotsc,p_r$ of the zeros of random polynomials as in \cite{ArmentanoBeltranShub2011} and $s=\floor{\sqrt{4r/(9J^2)}}$.} This construction yields a logarithmic energy of the form
\[
\mathbb{E}_{\textup{zeros}}\Biggl(\text{any $r$; $s=\floor[\Bigg]{\sqrt{\frac{4r}{9J^2}}}$}\Biggr)= \kappa n^2-\frac{1}{3}n\log n+C_{\mathrm{zeros}}n +  o(n),
\]
where
\begin{equation}\label{eq:czeros}
C_{\mathrm{zeros}}=-\frac13\log\frac{4}{9J^2}+\frac{2}{3}+\frac12\log 2=0.9190058\ldots.
\end{equation}
This proves \cref{th:main}. By comparison with \eqref{eq:minimo}, we observe that the logarithmic energy coincides with the theoretical minimum up to the $O(n)$ term. As pointed out in Theorem \ref{th:main}, the first few values of $n$ that can be realized through this construction are $n=2,3,8,10,12,21,24,27,30,33,36,52,56\dotso$. Choosing slightly different values for $s$ (more exactly, any $s$ of the form $\sqrt{4r/(9J^2)}+o(\sqrt{r})$), one can get other values of $n$ keeping the same asymptotic value of the energy, since the difference falls into the $o(n)$ term.

For the proof of \cref{prop:zeros}, we follow the approach in \cite[Section 2]{delatorre2022expectedenergyzeroselliptic}. Let $\gamma_r(u)=\gamma_{r,1}(u)+\gamma_{r,2}(u)$, where
\begin{align*}
	\gamma_{r,1}(u)&=\frac{\Bigl(1-\frac{ru^2}{(1+u^2)^{r}-1}\Bigr)^2(1+u^2)^{r-2}}{(1+u^2)^r-1},\\
	\gamma_{r,2}(u)&=\frac{\Bigl(1-\frac{ru^2(1+u^2)^{r-1}}{(1+u^2)^r-1}\Bigr)^2}{(1+u^2)^r-1}.
\end{align*}

\begin{lemma}\label{lem:gamma}
	Let $r\geq2$. Then,
	\[
	\gamma_{r,1}(u)\leq 
	\begin{cases}
		11 ru^2, &0<u<r^{-1/2},\\
		\frac{2}{(1+u^2)^2}, &r^{-1/2}\leq u<\infty,
	\end{cases}
	\quad
	\gamma_{r,2}(u)\leq 
	\begin{cases}
		3 ru^2, & 0<u<r^{-1/2},\\
		\frac{8r^2u^4}{(1+u^2)^{r+2}} & r^{-1/2}\leq u<\infty.
	\end{cases}
	\]
	Moreover,
	\[
	\gamma_r(u)\leq 
	\begin{cases}
		14 ru^2 & 0<u<r^{-1/2},\\
		\frac{34}{(1+u^2)^2} & r^{-1/2}\leq u<\infty.
	\end{cases}
	\]
\end{lemma}

\begin{proof}
	See \cref{sec:proof-lem:gamma}.
\end{proof}

\begin{lemma}\label{lem:constanteJ}
	The following limit exists:
	\[
	\lim_{r\to\infty}r^{3/2}\biggl(2\int_0^\infty u\log\biggl(1+\frac{u}{\sqrt{1+u^2}}\biggr)\gamma_r(u)\dd u-\frac{1}{2}\biggr),
	\]
	and its value is
	\begin{equation}\label{eq:J}
		J=\int_0^\infty\sqrt{t}\,
		\frac{2e^{2t}-4e^t-4te^{2t}+t^2e^t+t^2e^{2t}+4te^t+2}{(e^t-1)^3}\dd t= -0.5787893\ldots.
	\end{equation}
\end{lemma}

\begin{proof}
	\cref{sec:proof-lem:constanteJ}.
\end{proof}

\begin{proof}[Proof of \cref{prop:zeros}]
	Let
	\[
	I_r=\expectation[\Bigg]{}{\sum_{i\neq i'}\log\Biggl(1+\sqrt{\frac{1-\innerprod{p_i}{p_{i'}}}{2}}\Biggr)},
	\]
	where $p_1,\ldots,p_r$ are the images under the inverse stereographic projection of the zeros of elliptic random polynomials of degree $r$. In \cite[Section 2]{delatorre2022expectedenergyzeroselliptic} a method is given to compute expected values of arbitrary energies of this set of points. In the case of our energy, their approach yields:
	\begin{equation*}
		\frac{I_r}{r^2}=2\int_0^\infty u\log\biggl(1+\frac{u}{\sqrt{1+u^2}}\biggr)\gamma_r(u)\dd u.
	\end{equation*}
	From Lemma \ref{lem:gamma}, the integrand is bounded above by a function of the form $C/(1+u^2)$ for some constant $C>0$. Since this function is integrable over $[0,\infty)$, we may apply the dominated convergence theorem to obtain
	\begin{equation*}
		\lim_{r\to\infty}\frac{I_r}{r^2}=2\int_0^\infty u\log\biggl(1+\frac{u}{\sqrt{1+u^2}}\biggr)\lim_{r\to\infty}\gamma_r(u)\dd u,
	\end{equation*}
	provided this last limit exists. Indeed, we observe that
	\[
	\lim_{r\to\infty}\gamma_{r,1}(u)=\frac{1}{(1+u^2)^2},\quad \lim_{r\to\infty}\gamma_{r,2}(u)=0,
	\]
	and hence $\lim_{r\tendsto\infty}\gamma_r(u)=1/(1+u^2)^2$. Then, using the change of variables $t=u/\sqrt{1+u^2}$, we deduce that
	\begin{equation}\label{eq:auxn}
		\lim_{r\to\infty}\frac{I_r}{r^2}=2\int_0^\infty \frac{u}{(1+u^2)^2}\log\biggl(1+\frac{u}{\sqrt{1+u^2}}\biggr)\dd u
		=2\int_0^1 t\log(1+t)\dd t
		=\frac{1}{2}.
	\end{equation}
	Thus, we have $I_r=r^2/2+o(r^2)$. To compute the next-order term, define
	\[
	J_r=\frac{1}{\sqrt{r}}\biggl(I_r-\frac{r^2}{2}\biggr).
	\]
	Lemma \ref{lem:constanteJ} asserts that the limit $J=\lim_{r\to\infty}J_r$ exists, which implies that
	\begin{equation*}
		I_r=\frac{r^2}{2}+\sqrt{r}{J_r}=\frac{r^2}{2}+\sqrt{r}J+o(\sqrt{r}).
	\end{equation*}
	Finally, from Theorem \ref{th:logenergy} we conclude that
	\begin{align*}
		\mathbb{E}_{\textup{zeros}}&=-\frac{n^2}{2}\log2+\frac{n}2\log2-n\log s-s^2I_r\\
		&=-\frac{n^2}{2}\log2+\frac{n}2\log2-n\log s-s^2\biggl(\frac{r^2}{2}+\sqrt{r}J+o(\sqrt{r})\biggr)\\
		&=\kappa n^2\log2-n\log s+\frac{n}2\log2-s^2\sqrt{r}J+o\bigl(s^2\sqrt{r}\,\bigr)\\
		&=\kappa n^2\log2-n\log s+\frac{n}2\log2-\sqrt{n}s^{3/2}J+o\bigl(\sqrt{n}s^{3/2}\bigr).
	\end{align*}
	The proposition follows.
\end{proof}

\subsection{Point set \ref{item:EAP}: points coming from equal area partitions}

We now assume that $p_1,\ldots,p_r$ are chosen randomly, each of them in one region of an equal area partition of the sphere, as constructed in \cite{Leopardi}. More specifically, the sphere $\s$ is partitioned into $r$ open, pairwise disjoint subsets $B_1,\dotsc,B_r$, each of relative area $1/r$, and such that the union of their boundaries has zero measure on $\s$ and each $B_i$ is contained in some ball $\set{q\in\s\st \innerprod{p_i}{q}\geq 1-C_0/r}$, with $C_0$ some fixed constant. Then, $p_i$ is chosen randomly inside $B_i$ for $i=1,\dotsc,r$.

\begin{proposition}\label{prop:eap}
Assume that $p_1,\ldots,p_r\in\s$ are chosen randomly, each of them in one region of an equal area partition of the sphere as described above. Then, the expected logarithmic energy of the corresponding collection of special orthogonal matrices associated with $p_1,\ldots,p_r$ and a positive integer $s$ satisfies
\[
\mathbb{E}_{\textup{EAP}}\leq\kappa n^2-n\log s+\frac{n}2\log2+\textup{constant}\cdot \sqrt{n}s^{3/2}+O(s^2).
\]
\end{proposition}

It is immediate to see that, for fixed $n=rs$, the expected energy in this case is minimized when $s=\textup{constant}\cdot n^{1/3}$, or, equivalently, $s=\textup{constant}\cdot \sqrt{r}$. For instance, choosing this constant equal to $1$ and taking the integer part yields our third complete construction of $n$ points on $\so$, valid for a certain infinite sequence of values of $n$: \emph{the collection of special orthogonal matrices associated with points $p_1,\ldots,p_{r}\in\S^2$ chosen randomly, each of them in one region of an equal area partition of the sphere, and $s=\floor{\sqrt{r}}$}. The expected energy of this construction is
\[
\mathbb{E}_{\textup{EAP}}(\text{any $r$; $s=\lfloor\sqrt{r}\rfloor$})=\kappa n^2-\frac13n\log n+O(n).
\]
By comparison with \eqref{eq:minimo}, we observe that the logarithmic energy of this construction matches the theoretical minimum up to the second-order term. An explicit constant giving an upper bound for the term $O(n)$ can be estimated with some effort. The first few values of $n$ that can be realized through this construction are $n=2,3,8,10,12,14,16,27,\dotso$.

\begin{proof}[Proof of Proposition \ref{prop:eap}]
	Fix $i\in\set{1,\dotsc,r}$, and let $f\from\s\to\R$ be any integrable function. Denoting by $\expectation{p\in B}{}$ the expected value over a uniformly chosen $p\in B$, we have
	\begin{align*}
		\sum_{i'\neq i}\expectation{p_{i'}\in B_{i'}}{f(p_{i'})}&=\frac{r}{\vol(\s)}\int_{p\in\union_{i'\neq i}B_{i'}}f(p)\dd p\\
		&=\frac{r}{\vol(\s)}\int_{p\in\s}f(p)\dd p-\frac{1}{\vol(B_i)}\int_{p\in B_{i}}f(p)\dd p.
	\end{align*}
	In particular, 
	\begin{align*}
		\MoveEqLeft\sum_{i'\neq i}\expectation[\Bigg]{p_{i'}\in B_{i'}}{\log\biggl(1+\sqrt{\frac{1-\innerprod{p_i}{p_{i'}}}{2}}\biggr)}\\
		&=\frac{r}{\vol(\s)}\int_{p\in\s}\log\biggl(1+\sqrt{\frac{1-\innerprod{p_i}{p}}2}\biggr)\dd p\\
		&\quad-\frac{1}{\vol(B_i)}\int_{p\in B_i}\log\biggl(1+\sqrt{\frac{1-\innerprod{p_i}{p}}2}\biggr)\dd p.
	\end{align*}
	By \cref{prop:continuousesfera}, the first term equals $r/2$. For the second term, since
	\begin{equation*}
		B_i\containedeq\set{q\in\s\st\innerprod{p_i}{q}\geq 1-C_0/r},
	\end{equation*}
	we have
	\[
	\frac{1}{\vol(B_i)}\int_{p\in B_i}\log\biggl(1+\sqrt{\frac{1-\innerprod{p_i}{p}}{2}}\biggr)\leq
	 \log\biggl(1+\sqrt{\frac{C_0}{2r}}\biggr)=\frac{\sqrt{C_0}}{\sqrt{2r}}+O(1/r).
	\]
	Therefore,
	\[
	\sum_{i'\neq i}\expectation[\Bigg]{p_{i'}\in B_{i'}}{\log\biggl(1+\sqrt{\frac{1-\innerprod{p_i}{p_{i'}}}{2}}\biggr)}\geq
	\frac{r}2-\frac{\sqrt{C_0}}{\sqrt{2r}}+O(1/r).
	\]
	Since this is valid for all $i$, summing over $i$ gives
	\begin{equation*}
		\sum_{i'\neq i}\expectation[\Bigg]{p_{i'}\in B_{i'}}{\log\biggl(1+\sqrt{\frac{1-\innerprod{p_i}{p_{i'}}}{2}}\biggr)}\geq\frac{r^2}{2}-\frac{\sqrt{C_0r}}{\sqrt{2}}+O(1).
	\end{equation*}
	Finally, applying \cref{th:logenergy}, the expected logarithmic energy of the corresponding collection of special orthogonal matrices associated with $p_1,\dotsc,p_r$ satisfies
	\begin{align*}
		\mathbb{E}_{\textup{EAP}}&\leq-\frac{n^2}{2}\log2+\frac{n}2\log2-n\log s-s^2\biggl(\frac{r^2}{2}-\frac{\sqrt{C_0r}}{\sqrt{2}}+O(1)\biggr)\\
		&=-\frac{1+\log2}{2}n^2+\frac{n}2\log2-n\log s+\frac{\sqrt{C_0}}{\sqrt{2}}\sqrt{n}s^{3/2}+O(s^2),
	\end{align*}
	which completes the proof.
\end{proof}

\subsection{Point set \ref{item:sph-ensemble}: the spherical ensemble}

We now consider the case where the points $p_1,\dotsc,p_r$ are sampled from the spherical ensemble, as described in \cite{AlishahiZamani}. Specifically, two random $r\times r$ complex matrices $A$ and $B$ are generated, with entries drawn independently from the complex Gaussian distribution $\mathcal N_\mathbb{C}(0,1)$. The $r$ eigenvalues of the matrix $A^{-1}B$ are then computed and mapped to the sphere $\s$ using the inverse stereographic projection.

\begin{proposition}\label{prop:sphericalensemble}
Assume that $p_1,\dotsc,p_r\in\s$ are chosen from the spherical ensemble as in \cite{AlishahiZamani}. Then, the expected logarithmic energy of the corresponding collection of special orthogonal matrices associated with $p_1,\dotsc,p_r$ and a positive integer $s$ is
\begin{align*}
	\mathbb{E}_{\textup{sph}}&=\kappa n^2-n\log  s+\frac{\log2}{2}n+\frac{\sqrt\pi}{2}\frac{\Gamma(n/s)}{\Gamma(n/s+1/2)}s n - \frac{s^2}{2}\\
	&=\kappa n^2-n\log  s+\frac{\log2}{2}n+\frac{\sqrt  \pi s^{3/2}}{2}\sqrt n +  O(s^2)+o(s^{3/2}\sqrt n).
\end{align*}
\end{proposition}

It is straightforward to see that, for fixed $n=rs$, the minimum of the higher-order terms in this quantity is attained when $s=(16n/(9\pi))^{1/3}$ or, equivalently, $s=\sqrt{16r/(9\pi)}$. Since this last number is not an integer, we take its integer part to define our fourth complete construction of $n$ points on $\so$, valid for a certain infinite sequence of values of $n$: \emph{the collection of special orthogonal matrices associated with points $p_1,\dotsc,p_r$, sampled from the spherical ensemble, and $s=\floor{\sqrt{16r/(9\pi)}}$}. The expected logarithmic energy of this construction satisfies
\[
\mathbb{E}_{\textup{sph}}\biggl(\text{any $r$; $s=\floor[\bigg]{\sqrt{\frac{16r}{9\pi}}}$}\biggr)= \kappa n^2-\frac{1}{3}n\log n+C_{\textup{sph}}n +  o(n),
\]
where
\[
C_{\textup{sph}}=-\frac{5}{6}\log2+\frac{2}{3}\log3+\frac{1}{3}\log\pi+\frac{2}{3}=1.203028\ldots.
\]
Comparing with \eqref{eq:minimo}, we see that the logarithmic energy in this case matches the theoretical minimum up to the $O(n)$ term. The first few values of $n$ that can be realized through this construction are $n=2,3,4,5,6,7,16,18,20,\dotso$.

\begin{proof}[Proof of Proposition \ref{prop:sphericalensemble}]
	Following \cite[Eq. (4.4)]{AlishahiZamani}, for any measurable function $\phi\from\s\times\s\to[0,\infty)$, we have
	\[
	\expectation[\bigg]{}{\sum_{i\neq i'}\phi(p_i,p_{i'})}=\frac{r^2}{(4\pi)^2}\int_{p,q\in \s}\phi(p,q)\biggl(1-\biggl(1-\frac{\norm{p-q}^2}{4}\biggr)^{r-1}\biggr)\dd p\dd q.
	\]
	If $\phi(p,q)=f(\innerprod{p}{q})$ for some nonnegative, measurable function $f$, then by rotational symmetry we may fix $q=e_3$. Then, using that $\norm{p-e_3}^2=2-2p_3$, we obtain
	\[
	\expectation[\bigg]{}{\sum_{i\neq i'}\phi(p_i,p_{i'})}=\frac{r^2}{4\pi}\int_{p\in \s}f(p_3)\Biggl(1-\biggl(\frac{1+p_3}{2}\biggr)^{r-1}\Biggr)\dd p.
	\]
	Therefore,
	\begin{align*}
		\MoveEqLeft\expectation[\Bigg]{}{\sum_{i\neq i'}\log\biggl(1+\sqrt{\frac{1-\innerprod{p_i}{p_{i'}}}{2}}\biggr)}\\
		&=\frac{r^2}{4\pi}\int_{p\in \s}\log\biggl(1+\sqrt{\frac{1-p_3}{2}}\biggr)\Biggl(1-\biggl(\frac{1+p_3}{2}\biggr)^{r-1}\Biggr)\dd p.
	\end{align*}
	As in the proof of \cref{prop:continuousesfera}, we can reduce this last expresion to a one-dimensional integral over $t\in[-1,1]$. If, in addition, we use the change of variables $t=1-2u$ and integrate by parts, we get
	\begin{align*}
		\MoveEqLeft\expectation[\Bigg]{}{\sum_{i\neq i'}\log\biggl(1+\sqrt{\frac{1-\innerprod{p_i}{p_{i'}}}{2}}\biggr)}\\
		&=\frac{r^2}{2}\int_{-1}^1\log\biggl(1+\sqrt{\frac{1-t}{2}}\biggr)\Biggl(1-\biggl(\frac{1+t}{2}\biggr)^{r-1}\Biggr)\dd t\\
		&=r^2\int_{0}^1\log\bigl(1+\sqrt{u}\bigr)\Bigl(1-(1-u)^{r-1}\Bigr)\dd u\\
		&=r^2\biggl(\log2-\frac12\int_0^1\frac{u+\frac{1}{r}(1-u)^r}{u+\sqrt{u}}\dd u\biggr)\\
		&=r^2\biggl(\log2-\int_0^1\frac{x^2+\frac{1}{r}(1-x^2)^r}{1+x}\dd x\biggr)\\
		&=r^2\biggl(\frac12-\frac{1}{r}\int_0^1 (1-x)(1-x^2)^{r-1}\dd x\biggr)\\
		&=\frac{r^2}{2}-r\int_0^1 (1-x^2)^{r-1}\dd x+r\int_0^1 x(1-x^2)^{r-1}\dd x\\
		&=\frac{r^2}{2}-\frac{r}{2}\Beta\biggl(\frac12,r\biggr)+\frac{r}{2}\Beta(1,r)\\
		&=\frac{r^2}{2}-\frac{\sqrt \pi}{2}\frac{r\Gamma(r)}{\Gamma(r+1/2)}+\frac{1}{2},
	\end{align*}
	where $\Beta(x,y)$ is the beta function and we have used \cite[8.380-1]{integrales}.
	
	Using \cref{th:logenergy}, the expected logarithmic energy for the spherical ensemble on $\so$ is equal to
	\begin{align*}
		-\frac{n^2}{2}\log2+\frac{n}2\log2-n\log s-s^2\biggl(\frac{r^2}{2}-\frac{\sqrt \pi}{2}\frac{r\Gamma(r)}{\Gamma(r+1/2)}+\frac{1}{2}\biggr).
	\end{align*}
	The proposition follows from the previous equation after rearranging terms. The asymptotic expression in the proposition is a direct consecuence of Stirling's approximation for the Gamma function.
\end{proof}

\subsection{Point set \ref{item:harmonic}: the harmonic ensemble}

The harmonic ensemble introduced in \cite{BeltranMarzoOrtega} defines, for any positive integer $r=(L+1)^2$, which is a perfect square, a random collection of points $p_1,\ldots,p_r$ on the sphere with the property that 
\begin{align}
	\MoveEqLeft\expectation[\Bigg]{}{\sum_{i\neq i'}\log\biggl(1+\sqrt{\frac{1-\innerprod{p_i}{p_{i'}}}{2}}\biggr)}\nonumber\\
	&=\frac{1}{16\pi^2}\int_{p,p'\in\s}H_L(p,p')\log\biggl(1+\sqrt{\frac{1-\innerprod{p}{p'}}2}\biggr)\dd p\dd p'\label{eq:harmonicintegral},
\end{align}
see \cite[Proposition 1]{BeltranMarzoOrtega}. Here, $H_L(p,p')$ is defined in terms of the projection kernel $K_L(p,p')$ onto the subspace spanned by the spherical harmonics of degree at most $L$, which is a subspace of $L^2(\s)$ of dimension $r=(L+1)^2$, by the following formula:
\begin{equation}\label{eq:HLharmonic}
H_L(p,p')=K_L(p,p)K_L(p',p')-\abs{K_L(p,p')}^2.
\end{equation}
This is essentially the defining property of determinantal point processes; see, for example, \cite[Proposition 1]{BeltranMarzoOrtega} for examples of computations of different energies of spherical points using this approach, and \cite{HoughKrishnapurPeresVirag2009} for a comprehensive monograph on the topic.

\begin{proposition}\label{prop:harmonic}
	Let $p_1,\dotsc,p_r\in\s$ be $r$ points drawn from the harmonic ensemble on $\S^2$, with $r=(L+1)^2$. Then, the expected logarithmic energy of the corresponding collection of special orthogonal matrices associated with $p_1,\dotsc,p_r$ and a positive integer $s$ is given by
	\[
	\mathbb{E}_{\textup{$\s$-harmonic}}=\kappa n^2-n\log s+\frac{n}{2}\log2+g(n,s),
	\]
	where $g(n,s)$ is some function such that
	\begin{equation*}
		Cn^{1/2}s^{3/2}\log(n/s)\leq g(n,s)\leq C'n^{1/2}s^{3/2}\log(n/s),
	\end{equation*}
	for some positive constants $C<C'$.
\end{proposition}
	
It is immediate to see that, for fixed $n=rs$, the minimum of this quantity is attained when $s=o(n)$. This allows us to replace $\log(n/s)$ with $\log n$ in the bounds for $g(n,s)$.  Then, optimizing with respect to $s$ we find that a reasonable value for $s$ is %$s=n^{1/3}/(\log n)^{2/3}$, or, equivalently, 
$s=\sqrt{r}/\log r$. This yields our fifth complete construction of $n=r\floor{\sqrt{r}/\log r}$ points on $\so$: \emph{the collection of special orthogonal matrices associated with $r=(L+1)^2$ points $p_1,\dotsc,p_r$ drawn from the harmonic ensemble on $\S^2$ and $s=\floor{\sqrt{r}/\log r}$}. This construction satisfies $s=O(n^{1/3}/(\log n)^{2/3})$ and its logarithmic energy is
\[
\mathbb{E}_{\textup{$\s$-harmonic}}(\text{any $r$; $s=\floor{\sqrt{r}/\log r}$})=\kappa n^2-\frac13n\log n+O(n\log\log n).
\]
Comparing with \eqref{eq:minimo}, we observe that the logarithmic energy of this construction matches the theoretical minimum up to the second-order term, but it does not get to produce an $O(n)$ term. The first few values of $n$ that can be realized through this construction are $n=4,3,4,5,6,7,8\dotsc,72,73,74,150,152,154,156,\dotso$. 

To prove \cref{prop:harmonic}, we will use the following auxiliary result.

\begin{lemma}\label{lemma:integral-Jacobi}
	The following identity holds:
	\begin{equation*}
		\int_{0}^{\pi}\Jacobi{L}{1}{0}(\cos\theta)^2\log\bigl(1+\sin(\theta/2)\bigr)\sin\theta\dd \theta\asymp\frac{\log L}{L}.
	\end{equation*}
\end{lemma}

\begin{proof}
	See \cref{sec:proof-lemma:integral-Jacobi}.
\end{proof}

\begin{proof}[Proof of Proposition \ref{prop:harmonic}]
	It is known that the projection kernel $K_L$ can be expressed in terms of the Jacobi polynomial $\Jacobi{L}{1}{0}$ as follows (see \cite[Section 1.3]{BeltranMarzoOrtega}):
	\[
	K_L(p,p')=(L+1)\Jacobi{L}{1}{0}(\innerprod{p}{p'}).
	\]
	Hence, using $\Jacobi{L}{1}{0}(1)^2=(L+1)^2=r$,
	\begin{align*}
		H_L(p,p')&=(L+1)^2(\Jacobi{L}{1}{0}(1)^2-\Jacobi{L}{1}{0}(\innerprod{p}{p'})^2)\\
		&=r^2\biggl(1-\frac{1}{(L+1)^2}\Jacobi{L}{1}{0}(\innerprod{p}{p'})^2\biggr).\\
	\end{align*}
	Therefore, from \eqref{eq:harmonicintegral} and \eqref{eq:HLharmonic},
	\begin{align}
		\MoveEqLeft\expectation[\Bigg]{}{\sum_{i\neq i'}\log\biggl(1+\sqrt{\frac{1-\innerprod{p_i}{p_{i'}}}{2}}\biggr)}\label{eq:auxiliar}\\
			&=\frac{1}{16\pi^2}\int_{p,p'\in\s}r^2\biggl(1-\frac{1}{(L+1)^2}\Jacobi{L}{1}{0}(\innerprod{p}{p'})^2\biggr)\log\biggl(1+\sqrt{\frac{1-\innerprod{p}{p'}}2}\biggr)\dd p\dd p'\notag\\
			&=\frac{1}{4\pi}\int_{p\in\s}r^2\biggl(1-\frac{1}{(L+1)^2}\Jacobi{L}{1}{0}(p_3)^2\biggr)\log\biggl(1+\sqrt{\frac{1-p_3}2}\biggr)\dd p\notag\\
		&=\frac{r^2}{2}\int_{-1}^1\biggl(1-\frac{1}{(L+1)^2}\Jacobi{L}{1}{0}(t)^2\biggr)\log\biggl(1+\sqrt{\frac{1-t}2}\biggr)\dd t,\notag
	\end{align}
	where for the last equality we have again translated the integral in $\s$ into a univariate integral. Using the computations in the proof of \cref{prop:continuousesfera}, we can rewrite the previous expression as
	\[
	\frac{r^2}{2}-\frac{r}{2}\int_{-1}^1\Jacobi{L}{1}{0}(t)^2\log\biggl(1+\sqrt{\frac{1-t}2}\biggr)\dd t.
	\]
	Applying the change of variables $t=\cos\theta$, the integral becomes
	\[
	\int_{0}^{\pi}\Jacobi{L}{1}{0}(\cos\theta)^2\log\bigl(1+\sin(\theta/2)\bigr)\sin\theta\dd \theta.
	\]
	From \cref{lemma:integral-Jacobi}, there exist positive constants $C$ and $C'$ such that
	\[
	\frac{r^2}{2}-C\sqrt{r}\log r\leq \expectation[\Bigg]{}{\sum_{i\neq i'}\log\biggl(1+\sqrt{\frac{1-\innerprod{p_i}{p_{i'}}}{2}}\biggr)}\leq\frac{r^2}{2}-C'\sqrt{r}\log r.
	\]
	Finally, from Theorem \ref{th:logenergy} and using that $n=rs$, we deduce that
	\begin{align*}
		\mathbb{E}_{\textup{$\s$-harmonic}}&\leq-\frac{n^2}{2}\log2+\frac{n}2\log2-n\log s-s^2\biggl(\frac{r^2}{2}-C\sqrt{r}\log r\biggr)\\
		&=\kappa n^2-n\log s+\frac{n}{2}\log 2+C\sqrt{n}s^{3/2}\log\frac{n}{s},
	\end{align*}
	with the reverse inequality (replacing $C$ by $C'$) for the lower bound. The proposition follows.
\end{proof}

\section{Proofs of some auxiliary results}\label{sec:auxiliary-results}

\subsection{Auxiliary results for the proof of \cref{thm:energy-harmonic}}

\subsubsection{Proof of \cref{lemma:aux-integral-1}}\label{sec:proof-lemma:aux-integral-1}

To prove this lemma, we use Hilb's formula \eqref{eq:Hilb-Gegenbauer-2} for the Gegenbauer polynomials $\Gegenbauer{2L}{2}$. Using this formula, we have
\begin{equation*}
	\Gegenbauer{2L}{2}(\cos\theta)^2\sin^2\theta=A_1+A_2+A_3,
\end{equation*}
where
\begin{align*}
	A_1&=\frac{\pi}{8}(2L)^{3}(1+O(L^{-1}))\theta^{-1}\BesselJ{3/2}\bigl((2L+2)\theta\bigr)^2,\\
	A_2&=\textup{Error}^2,\\
	A_3&=L^{3/2}(1+O(L^{-1}))\theta^{-1/2}\BesselJ{3/2}\bigl((2L+2)\theta\bigr)\textup{Error},\\
\end{align*}
and
\begin{equation*}
	\textup{Error}=\begin{cases}
		O(L^{-3/2}), & L^{-1}\leq \theta\leq L^{-1/2},\\
		(\theta^{3})O(L^{3/2}), & 0<\theta<L^{-1}.
	\end{cases}
\end{equation*}
Hence, calling $I$ the integral of the lemma, we have
\begin{equation*}
	I=I_1+I_2+I_3,
\end{equation*}
where
\begin{equation*}
	I_j=\int_{0}^{L^{-1/2}}A_j\dd\theta.
\end{equation*}
First, note that $\textup{Error}$ is always $O(L^{-3/2})$. Hence, $I_2=O(L^{-7/2})$. To estimate $I_3$, note that $\abs{\BesselJ{3/2}(x)}\leq 1$ for all $x\in\R$, so
\begin{equation*}
	I_3=O\biggl(\int_{0}^{L^{-1/2}}\theta^{-1/2}\dd\theta\biggr)=O\biggl(\int_{0}^{L^{-1/2}}\theta^{-1/2}\dd\theta\biggr)=O(L^{-1/4}).
\end{equation*}
Finally, we estimate $I_1$. Denoting $\hat A_1=	\pi L^3\theta^{-1}\BesselJ{3/2}\bigl((2L+2)\theta\bigr)^2$ we have
\begin{equation*}
	I_1=(1+O(L^{-1}))\biggl(\int_{0}^{\infty}\hat A_1\dd\theta-\int_{L^{-1/2}}^{\infty}\hat A_1\dd\theta\biggr).
\end{equation*}
Using that $\BesselJ{3/2}(x)=O(x^{-1/2})$ as $x\tendsto\infty$,
\begin{equation*}
	\int_{L^{-1/2}}^{\infty}\hat A_1\dd\theta=O\biggl(L^3(2L+2)^{-1}\int_{L^{-1/2}}^{\infty}\theta^{-2}\dd\theta\biggr)=O(L^{5/2}).
\end{equation*}
Finally, we get
\begin{equation*}
	\int_{0}^{\infty}\hat A_1\dd\theta=\pi L^3\int_{0}^{\infty}\frac{\BesselJ{3/2}((2L+2)\theta)^2}{\theta}\dd\theta=\frac{\pi}{3}L^3.
\end{equation*}
where in the last integral involving the Bessel function we have used \eqref{eq:integral-Bessel}. The lemma follows. \qed

\subsubsection{Proof of \cref{lemma:aux-integral-2}}\label{sec:proof-lemma:aux-integral-2}

From \eqref{eq:bounds-Gegenbauer} we have the estimate $\Gegenbauer{2L}{2}=
\theta^{-2}O(L)$. Hence, using that $\sin\theta=O(\theta)$,
\begin{equation*}
	\abs[\bigg]{\int_{L^{-1/2}}^{\pi/2}\Gegenbauer{2L}{2}(\cos\theta)^2\sin^2\theta\dd\theta}=O\biggl(L^2\int_{L^{-1/2}}^{\pi/2}\theta^{-2}\biggr)=O(L^{5/2}).\pushQED{\qed}\qedhere
\end{equation*}

\subsubsection{Proof of \cref{lemma:aux-integral-3}}\label{sec:proof-lemma:aux-integral-3}

We use the estimates $\Gegenbauer{2L}{2}=
\theta^{-2}O(L)$ and $\sin\theta=O(\theta)$, together with $\abs{\log\sin\theta}\leq2\abs{\log \theta}$ for small $\theta$. Hence, the integral of the lemma can be estimated as
\begin{align*}
	\MoveEqLeft O\biggl(L^2\int_{L^{-1/2}}^{\pi/2}\theta^{-2}\log\theta\dd\theta\biggr)=O\biggl(L^2\biggl(-\frac{1}{\theta}-\frac{\log \theta}{\theta}\biggr)\eval_{L^{-1/2}}^{\pi/2}\biggr)=O(L^{5/2}\log L).\pushQED{\qed}\qedhere
\end{align*}

\subsubsection{Proof of \cref{lemma:aux-integral-4}}\label{sec:proof-lemma:aux-integral-4}

Note that in the range $[0,L^{-1/2}]$ we have
\begin{equation*}
	\log(\sin\theta)=\log(\theta(1+O(\theta^2)))=\log\theta+\log(1+O(\theta^2))=\log\theta+O(L^{-1}).
\end{equation*}
From Lemma \ref{lemma:aux-integral-1}, the integral of the lemma equals
\begin{equation*}
	\int_{0}^{L^{-1/2}}\log\theta\Gegenbauer{2L}{2}(\cos\theta)^2\sin^2\theta\dd\theta+O(L^2).
\end{equation*}
Let us estimate the integral in the previous expression. Making the substitution $\theta=\hat\theta/(2L)$, this integral equals
\begin{align*}
	\MoveEqLeft\frac{1}{2L}\int_{0}^{2L^{1/2}}\log(\hat\theta/(2L))\Gegenbauer{2L}{2}(\cos(\hat\theta/(2L)))^2\sin^2(\hat\theta/(2L))\dd\hat\theta\\
	&=\frac{1}{2L}\int_{0}^{2L^{1/2}}\log\hat\theta\Gegenbauer{2L}{2}(\cos(\hat\theta/(2L)))^2\sin^2(\hat\theta/(2L))\dd\hat\theta\\
	&\quad-\frac{\log (2L)}{2L}\int_{0}^{2L^{1/2}}\Gegenbauer{2L}{2}(\cos(\hat\theta/(2L)))^2\sin^2(\hat\theta/(2L))\dd\hat\theta\\
	&=\frac{1}{2L}\int_{0}^{2L^{1/2}}\log\hat\theta\Gegenbauer{2L}{2}(\cos(\hat\theta/(2L)))^2\sin^2(\hat\theta/(2L))\dd\hat\theta\\
	&\quad-\log (2L)\int_{0}^{L^{-1/2}}\Gegenbauer{2L}{2}(\cos\theta)^2\sin^2\theta\dd\theta.
\end{align*}
From \cref{lemma:aux-integral-1}, we have
\begin{align*}
	-\log (2L)\int_{0}^{L^{-1/2}}\Gegenbauer{2L}{2}(\cos\theta)^2\sin^2\theta\dd\theta&=-\log(2L)\biggl(\frac{\pi}{3}L^3+O(L^{5/2})\biggr)\\
	&=-\frac{\pi}{3}L^3\log (2L)+O(L^{5/2}\log L)\\
	&=-\frac{\pi}{3}L^3\log L-\frac{\pi\log2}{3}L^3+o(L^3).
\end{align*}
Finally, we claim that
\begin{align*}
	\MoveEqLeft\frac{1}{2L}\int_{0}^{2L^{1/2}}\log\hat\theta\Gegenbauer{2L}{2}(\cos(\hat\theta/(2L)))^2\sin^2(\hat\theta/(2L))\dd\hat\theta\\
	&=\frac{\pi}{8}\biggl(\frac{7-3\gamma-3\log{2}}{9}\biggr)(2L)^3+o(L^3).
\end{align*}
To prove this claim, it suffices to check that
\begin{equation*}
	\lim_{L\tendsto\infty}\frac{1}{(2L)^4}\int_{0}^{2L^{1/2}}\log\hat\theta\Gegenbauer{2L}{2}(\cos(\hat\theta/(2L)))^2\sin^2(\hat\theta/(2L))\dd\hat\theta=\frac{\pi}8\cdot\frac{7-3\gamma-3\log{2}}{9}
\end{equation*}
Note that we can write the previous expression as
\begin{align*}
	\MoveEqLeft\lim_{L\tendsto\infty}\frac{1}{(2L)^6}\int_{0}^{\infty}\chi_{[0,2L^{1/2}]}\log\hat\theta\Gegenbauer{2L}{2}(\cos(\hat\theta/(2L)))^2\hat\theta^2(1+o(1))\dd\hat\theta.
\end{align*}
Using the dominated convergence theorem and the Mehler--Heine formula \eqref{eq:Mehler-Heine-Gegenbauer-2}, this limit equals
\begin{equation*}
	\frac{\pi}{8}\int_{0}^{\infty}\frac{\BesselJ{3/2}(\hat\theta)^2\log\hat\theta}{\hat\theta}\dd\hat\theta=\frac{\pi}{8}\cdot\frac{7-3\gamma-3\log{2}}{9},
\end{equation*}
where the last equality follows from \cref{eq:integral-Bessel-log}.\qed

\subsection{Auxiliary results for the proof of \cref{prop:zeros}}

\subsubsection{Proof of \cref{lem:gamma}}\label{sec:proof-lem:gamma}

We apply Taylor's theorem with the Lagrange form of the remainder to the function $g(t)=(1+t)^r$, yielding, for some $\zeta\in(0,t)$,
\[
(1+t)^r=1+rt+\frac{r(r-1)(1+\zeta)^{r-2}}{2}t^2.
\]
Applying this with $t=u^2$, we obtain the following bounds for $0<u<r^{-1/2}$:
\begin{align*}
	(1+u^2)^r&\leq1+ru^2+\frac{r(r-1)(1+1/r)^{r-2}}{2}u^4\leq 1+ru^2+2r^2u^4,\\
	(1+u^2)^r&\geq 1+ru^2+\frac{r(r-1)}2u^4\geq 1+ru^2+\frac{r^2}{4}u^4,
\end{align*}
From these bounds, we have
\begin{align*}
	\gamma_{r,1}(u)&\leq \frac{\Bigl(1-\frac{ru^2}{ru^2+2r^2u^4}\Bigr)^2e}{ru^2}\leq \frac{4eru^2}{(1+2ru^2)^2}\leq 11ru^2,\\
	\gamma_{r,2}(u)&\leq \frac{\Bigl(1-\frac{ru^2(1+ru^2-u^2)}{ru^2+2r^2u^4}\Bigr)^2}{ru^2} =
	\frac{u^2(r+1)^2}{r(2ru^2+1)^2}\leq 3ru^2.
\end{align*}
Now consider the range $r^{-1/2}<u<\infty$. Since the function $x\mapsto x/(x-1)$ is decreasing and $(1+u^2)^r\geq (1+1/r)^r\geq 2$, we have
\[
\frac{(1+u^2)^r}{(1+u^2)^r-1}\leq 2,
\]
which implies that
\begin{equation*}
	\gamma_{r,1}(u)\leq\frac{1}{(1+u^2)^2}\frac{(1+u^2)^r}{(1+u^2)^r-1}\leq\frac{2}{(1+u^2)^2},
\end{equation*}
and
\begin{align*}
	\gamma_{r,2}(u)&\leq\frac{r^2u^4}{(1+u^2)^2((1+u^2)^r-1)}\biggl(\frac{(1+u^2)^r}{(1+u^2)^r-1}\biggr)^2\\
	&\leq\frac{4r^2u^4}{(1+u^2)^2((1+u^2)^r-1)}\\
	&\leq \frac{8r^2u^4}{(1+u^2)^2(1+u^2)^r}.
\end{align*}
The lemma follows. For the last bound on $\gamma$, recall that $\gamma=\gamma_{r,1}+\gamma_{r,2}$ and note that
\begin{align*}
	\frac{8r^2u^4}{(1+u^2)^{r+2}}&\leq \frac{8r^2u^4}{(1+u^2)^2(1+ru^2+r(r-1)u^4/2)}\\
	&\leq \frac{8r^2u^4}{(1+u^2)^2r^2u^4/4}\\
	&=\frac{32}{(1+u^2)^2}.\pushQED{\qed}\qedhere
\end{align*}

\subsubsection{Proof of \cref{lem:constanteJ}}\label{sec:proof-lem:constanteJ}

Let $J$ be the limit we want to compute. From the integral computed in \eqref{eq:auxn}, we may decompose
\begin{equation*}
	J=J_1+J_2,
\end{equation*}
where
\begin{align*}
	J_1&=\lim_{r\to\infty}\biggl(2\int_0^\infty r^{3/2}u\log\biggl(1+\frac{u}{\sqrt{1+u^2}}\biggr)\biggl(\gamma_{r,1}(u)-\frac{1}{(1+u^2)^2}\biggr)\dd u\biggr),\\
	J_2&=\lim_{r\to\infty}\biggl(2\int_0^\infty r^{3/2}u\log\biggl(1+\frac{u}{\sqrt{1+u^2}}\biggr)\gamma_{r,2}(u)\dd u\biggr).
\end{align*}
We begin with $J_2$. Let
\begin{align*}
	A_r&=\int_0^{1} r^{3/2}u\log\biggl(1+\frac{u}{\sqrt{1+u^2}}\biggr)\gamma_{r,2}(u)\dd u,\\
	B_r&=\int_1^\infty r^{3/2}u\log\biggl(1+\frac{u}{\sqrt{1+u^2}}\biggr)\gamma_{r,2}(u)\dd u.
\end{align*}
Making the substitution $s=u\sqrt r$, we rewrite $A_r$ as
\begin{equation*}
	A_r=\int_0^{\sqrt{r}} r^{1/2}s\log\biggl(1+\frac{s/\sqrt{r}}{\sqrt{1+s^2/r}}\biggr)\gamma_{r,2}(s/\sqrt{r})\dd s.
\end{equation*}
From Lemma \ref{lem:gamma}, the integrand is bounded above by a constant independent of $r$ in the interval $s\in[0,1]$, and in the interval $[1,\infty)$ it is bounded above by
\[
\frac{8s^6}{(1+s^2/r)^{r+2}}\leq\frac{8s^6}{(1+s^2/r)^{r}}\leq  \frac{8s^6}{(1+s^2/4)^{4}},
\]
(for $r\geq4$) which is integrable over $[1,\infty)$. Therefore, using the dominated convergence theorem we obtain
\begin{align*}
	\lim_{r\to\infty}A_r&=\int_0^\infty\lim_{r\to\infty}\biggl(r^{1/2}s\log\biggl(1+\frac{s/\sqrt{r}}{\sqrt{1+s^2/r}}\biggr)\gamma_{r,2}(s/\sqrt{r})\biggr)\dd s\\
	&=\int_0^\infty s^2\frac{\bigl(e^{s^2}-1-s^2e^{s^2}\bigr)^2}{\bigl(e^{s^2}-1\bigr)^3}\dd s\\
	&=\frac12\int_0^\infty\sqrt{t}\,\frac{(e^{t}-1-te^{t})^2}{(e^{t}-1)^3}\dd t.
\end{align*}
It is much easier to check that $B_r\tendsto0$. Indeed, from \cref{lem:gamma}, the integrand in $B_r$ is bounded above by
\[
\frac{8r^{7/2}u}{(1+u^2)^r}= \frac{8r^{7/2}}{(1+u^2)^{r-2}}\frac{u}{(1+u^2)^2}\leq\frac{C}{u^3},
\]
for some constant $C>0$ and all $u\in[1,\infty)$. Since $C/u^3$ is integrable over $[1,\infty)$, we may apply again the dominated convergence theorem, which yields
\begin{equation*}
	\lim_{r\to\infty}B_r=\int_1^\infty u\log\biggl(1+\frac{u}{\sqrt{1+u^2}}\biggr)\lim_{r\to\infty}\bigl(r^{3/2}\gamma_{r,2}(u)\bigr)\dd u=0.
\end{equation*}
All together, we have proved that
\[
J_2=\int_0^\infty\sqrt{t}\,\frac{(e^{t}-1-te^{t})^2}{(e^{t}-1)^3}\dd t.
\]
We now turn to the analysis of $J_1$. The integral in the definition of $J_1$ equals
\[
\int_0^\infty \frac{r^{3/2}u}{(1+u^2)^2}\log\biggl(1+\frac{u}{\sqrt{1+u^2}}\biggr)\left(\frac{\Bigl(1-\frac{ru^2}{(1+u^2)^{r}-1}\Bigr)^2(1+u^2)^{r}}{(1+u^2)^r-1}-1\right)\dd u.
\]
After some simple algebraic manipulations, we can rewrite this expression as
\[
\int_0^\infty \frac{r^{3/2}u}{(1+u^2)^2}\log\biggl(1+\frac{u}{\sqrt{1+u^2}}\biggr)\frac{1+(1+u^2)^{r}\Bigl(\frac{r^2u^4}{((1+u^2)^{r}-1)^2} -\frac{2ru^2}{(1+u^2)^{r}-1}\Bigr)}{(1+u^2)^r-1}\dd u.
\]
As before, we split this into two parts:
\begin{align*}
	C_r&=\int_0^1 \frac{r^{3/2}u}{(1+u^2)^2}\log\biggl(1+\frac{u}{\sqrt{1+u^2}}\biggr)\frac{1+(1+u^2)^{r}\Bigl(\frac{r^2u^4}{((1+u^2)^{r}-1)^2} -\frac{2ru^2}{(1+u^2)^{r}-1}\Bigr)}{(1+u^2)^r-1}\dd u,\\
	D_r&=\int_1^\infty \frac{r^{3/2}u}{(1+u^2)^2}\log\biggl(1+\frac{u}{\sqrt{1+u^2}}\biggr)\frac{1+(1+u^2)^{r}\Bigl(\frac{r^2u^4}{((1+u^2)^{r}-1)^2} -\frac{2ru^2}{(1+u^2)^{r}-1}\Bigr)}{(1+u^2)^r-1}\dd u.
\end{align*}
It is easy to see that $D_r\tendsto0$ as $r\tendsto\infty$. Indeed, note that the integrand in $D_r$ is bounded above by
\[
\textup{constant}\cdot \frac{r^{3/2}u}{(1+u^2)^2}\cdot \frac{2ru^2}{(1+u^2)^{r}}\leq \textup{constant}\cdot \frac{r^{5/2}u^3}{(1+u^2)^{2+r}}\leq\frac{\textup{constant}}{u^2},
\]
where the constant can be chosen independent of $r$. Hence, by the dominated convergence theorem we can interchange limit and integral and we get
\begin{equation*}
	0\leq \lim_{r\to\infty}D_r\leq \int_1^\infty \lim_{r\to\infty}\textup{constant}\cdot \frac{r^{5/2}u^3}{(1+u^2)^{2+r}}\dd u=0.
\end{equation*}
It only remains to deal with $C_r$. As before, making the change of variables $s=u\sqrt r$ we obtain
\begin{align*}
	\MoveEqLeft C_r=\int_0^{\sqrt{r}} \frac{r^{1/2}s}{(1+s^2/r)^2}\log\biggl(1+\frac{s/\sqrt{r}}{\sqrt{1+s^2/r}}\biggr)\\
	&\qquad\times\frac{1+(1+s^2/r)^{r}\Bigl(\frac{s^4}{((1+s^2/r)^{r}-1)^2} -\frac{2s^2}{(1+s^2/r)^{r}-1}\Bigr)}{(1+s^2/r)^r-1}\dd s.
\end{align*}
Elementary manipulations show that the integrand is bounded above in the interval $[0,1]$ by
\[
s^2\frac{1+e\Bigl(\frac{s^4}{((1+s^2/2)^{2}-1)^2} -\frac{2s^2}{e^{s^2}-1}\Bigr)}{(1+s^2/2)^2-1}\leq \textup{constant}.
\]
For $s\in [1,\infty]$, using that $(1+s^2/r)^r-1\geq (1+s^2/r)^r/4$, we can bound it above by
\[
\textup{constant}\cdot \frac{s^2}{(1+s^2/r)^{r}}\biggl(\frac{s^4}{(1+s^2/r)^{r}} -2s^2\biggr)\leq \frac{\textup{constant}\cdot s^6}{(1+s^2/4)^4},
\]
when $r\geq4$. Putting everything together, we have bounded the integrand in $C_r$ by a function that is integrable over $[0,\infty)$, so we may again apply the dominated convergence theorem. This yields
\begin{align*}
	\lim_{r\to\infty}C_r&=\int_0^\infty s^2\frac{1+e^{s^2}\Bigl(\frac{s^4}{(e^{s^2}-1)^2}-\frac{2s^2}{e^{s^2}-1}\Bigr)}{e^{s^2}-1}\dd s\\
	&=\frac12\int_0^\infty \sqrt{t}\,\frac{1+e^{t}\Bigl(\frac{t^2}{(e^{t}-1)^2}-\frac{2t}{e^{t}-1}\Bigr)}{e^{t}-1}\dd t.
\end{align*}
In summary, we have proved that
\begin{align*}
	J&=J_1+J_2\\
	&=\lim_{r\to\infty}(2A_r+2B_r+2C_r+2D_r)\\
	&=\int_0^\infty\sqrt{t}\,\frac{(e^{t}-1-te^{t})^2}{(e^{t}-1)^3}\,dt + \int_0^\infty \sqrt{t}\frac{1+e^{t}\Bigl(\frac{t^2}{(e^{t}-1)^2}-\frac{2t}{e^{t}-1}\Bigr)}{e^{t}-1}\dd t\\
	&=\int_0^\infty\sqrt{t}\,
	\frac{2e^{2t}-4e^t-4te^{2t}+t^2e^t+t^2e^{2t}+4te^t+2}{(e^t-1)^3}\dd t.\pushQED{\qed}\qedhere
\end{align*}

\subsection{Auxiliary results for the proof of \cref{prop:harmonic}}

\subsubsection{Proof of \cref{lemma:integral-Jacobi}}\label{sec:proof-lemma:integral-Jacobi}

Let
\begin{equation*}
	I=\int_{0}^{\pi}\Jacobi{L}{1}{0}(\cos\theta)^2\log\bigl(1+\sin(\theta/2)\bigr)\sin\theta\dd \theta
\end{equation*}
be the integral in the lemma. First, we split this integral as
\begin{equation*}
	I=I_1+I_2+I_3,
\end{equation*}
where
\begin{align*}
	I_1&=\int_{0}^{L^{-1/2}}\Jacobi{L}{1}{0}(\cos\theta)^2\log\bigl(1+\sin(\theta/2)\bigr)\sin\theta\dd \theta,\\
	I_2&=\int_{L^{-1/2}}^{\pi/2}\Jacobi{L}{1}{0}(\cos\theta)^2\log\bigl(1+\sin(\theta/2)\bigr)\sin\theta\dd \theta,\\
	I_3&=\int_{\pi/2}^{\pi}\Jacobi{L}{1}{0}(\cos\theta)^2\log\bigl(1+\sin(\theta/2)\bigr)\sin\theta\dd \theta.
\end{align*}
For the integral $I_1$, using that $\sin\theta=O(\theta)$ and $\log\bigl(1+\sin(\theta/2)\bigr)=O(\theta)$ together with the estimates \eqref{eq:bounds-Jacobi}, we have
\begin{align*}
	I_1&=O\biggl(\int_0^{L^{-1}}L^2\theta^2\dd\theta+L^{-1}\int_{L^{-1}}^{L^{-1/2}}\theta^{-1}\dd\theta\biggr)\\
	&=O(L^{-1}\log L).
\end{align*}
To estimate $I_3$, making the substitution $\theta\mapsto\pi-\theta$ and using that
\begin{equation*}
	\Jacobi{L}{\alpha}{\beta}(-x)=(-1)^{L}\Jacobi{L}{\beta}{\alpha}(x),
\end{equation*}
the integral reads as
\begin{equation*}
	\int_{0}^{\pi/2}\Jacobi{L}{0}{1}(\cos\theta)^2\log\bigl(1+\cos(\theta/2)\bigr)\sin\theta\dd \theta.
\end{equation*}
Using that $\sin\theta=O(\theta)$ and $\log\bigl(1+\cos(\theta/2)\bigr)=O(1)$, together with the bounds \eqref{eq:bounds-Jacobi} for the Jacobi polynomials, we have
\begin{align*}
	I_3&=O\biggl(\int_0^{L^{-1}}\theta\dd\theta+L^{-1}\int_{L^{-1}}^{\pi/2}\theta^{-1}\theta\dd\theta\biggr)\\
	&=O(L^{-1}).\qedhere
\end{align*}
Finally, we estimate $I_2$. Note that
\begin{equation*}
	\log\bigl(1+\sin(\theta/2)\bigr)\sin\theta\asymp \theta^2.
\end{equation*}
Using this estimate and the bounds \eqref{eq:bounds-Jacobi} for the Jacobi polynomials, we have
\begin{equation*}
	I_2\lesssim L^{-1}\int_{L^{-1/2}}^{\pi/2}\theta^{-1}\dd\theta\lesssim \frac{\log L}{L}.
\end{equation*}
For the lower bound, recall that  we have the following asymptotic representation (see \cite[Formula 8.965]{integrales}):
\begin{equation*}
	\Jacobi{L}{1}{0}(\cos\theta)=\frac{\cos((L+1)\theta-3\pi/4)}{\sqrt{\pi L}\sin^{3/2}(\theta/2)\cos^{1/2}(\theta/2)}+O(L^{-3/2}).
\end{equation*}
Then, we have
\begin{equation*}
	\Jacobi{L}{1}{0}(\cos\theta)^2\gtrsim \frac{\cos^2((L+1)\theta-3\pi/4)}{\pi L\sin^{3}(\theta/2)\cos(\theta/2)}.
\end{equation*}
Note that, for $\theta\in[L^{-1/2},\pi/2]$,
\begin{equation*}
	\cos^2((L+1)\theta-3\pi/4)=1\iff \theta=\frac{k\pi+3\pi/4}{L+1},\qquad \floor{c\sqrt{L}}\leq k\leq\floor{c'L},
\end{equation*}
for certain constants $c,c'>0$. Hence,
\begin{equation*}
	\frac{k\pi+3\pi/4- \pi/4}{L+1}\leq \theta\leq \frac{k\pi+3\pi/4+ \pi/4}{L+1}\implies \cos^2\bigl((L+1)\theta-3\pi/4\bigr)\geq 1/2.
\end{equation*}
Then,
\begin{align*}
	I_2
	&\gtrsim\frac{1}{L} \sum_{k=\floor{c\sqrt{L}}}^{\floor{c'L}}\int_{\frac{k\pi+\pi/2}{L+1}}^{{\frac{k\pi+\pi}{L+1}}}\theta^{-1}\dd\theta=\frac{1}{L}\sum_{k=\floor{c\sqrt{L}}}^{\floor{c'L}}\log\biggl(1+\frac{\pi/2}{k\pi+\pi/2}\biggr)\\
	&\gtrsim\frac{1}{L}\sum_{k=\floor{c\sqrt{L}}}^{\floor{c'L}}\frac{1}{2k+1}\gtrsim \frac{1}{L}\int_{\floor{c\sqrt{L}}}^{\floor{c'L}+1}\frac{1}{2x+1}\dd x\\
	&=\frac{1}{2L}\log\biggl(\frac{2\floor{c'L}+3}{2\floor{c\sqrt{L}}+1}\biggr)\asymp \frac{\log(L)}{L}.\pushQED{\qed}\qedhere
\end{align*}

%----------- APPENDICES

\appendix

\section{Orthogonal polynomials and special functions}\label{appendix:orthogonal polynomials}

We compile in this appendix some properties about orthogonal polynomials that are used throughout this document.

\subsection{Jacobi polynomials}

Jacobi polynomials $\Jacobi{L}{\alpha}{\beta}$, where $L$ is the degree of the polynomial and $\alpha$ and $\beta$ are parameters, form a complete orthogonal system in $[-1,1]$ with respect to the weight $w(x)=(1-x)^{\alpha}(1+x)^{\beta}$.

From \cite[Eq. (7.32.5)]{Szego1975} we have the following asymptotic estimates for the Jacobi polynomials (for $\alpha\geq -1/2$):
\begin{equation}\label{eq:bounds-Jacobi}
	\Jacobi{L}{\alpha}{\beta}(\cos\theta)=\begin{cases}
		\theta^{-\alpha-1/2}O(L^{-1/2}), & L^{-1}\leq \theta\leq \pi/2, \\
		O(L^{\alpha}), & 0\leq \theta\leq L^{-1}.
	\end{cases}
\end{equation}

A more precise result to study the asymptotic behavior of Jacobi polynomials is the following formula of Hilb's type.

\begin{theorem}[{\cite[Theorem 8.21.12]{Szego1975}}]\label{thm:Hilb-Jacobi}
	Let $\alpha>-1$, and let $\beta$ be arbitrary and real. Then, we have
	\begin{equation*}
		\biggl(\sin\frac{\theta}{2}\biggr)^{\alpha}\biggl(\cos\frac{\theta}{2}\biggr)^{\beta}\Jacobi{L}{\alpha}{\beta}(\cos\theta)=A(L)^{-\alpha}\frac{\Gamma(L+\alpha+1)}{L!}\frac{\theta^{1/2}}{\sin^{1/2}\theta}\BesselJ{\alpha}\bigl(A(L)\theta\bigr)+\textup{Error},
	\end{equation*}
	where $\BesselJ{\alpha}$ is the Bessel function of the first kind of order $\alpha$, $A(L)=\bigl(L+(\alpha+\beta+1)/2\bigr)$, and
	\begin{equation*}
		\textup{Error}=\begin{cases}
			\theta^{1/2}O(L^{-3/2}), & L^{-1}\leq \theta\leq \pi-\varepsilon,\\
			\theta^{\alpha+2}O(L^{\alpha}), & 0<\theta<L^{-1},
		\end{cases}
	\end{equation*}
	where $\varepsilon$ is a fixed positive constant.
\end{theorem}

Near the endpoints, the asymptotic behavior of the Jacobi polynomials is given by the Mehler--Heine formula (see \cite[Theorem 8.1.1]{Szego1975}):
\begin{equation}\label{eq:Mehler-Heine}
	\lim_{L\tendsto\infty}L^{-\alpha}\Jacobi{L}{\alpha}{\beta}\biggl(\cos\frac{z}{L}\biggr)=(z/2)^{-\alpha}\BesselJ{\alpha}(z).
\end{equation}

\subsection{Gegenbauer polynomials}

Gegenbauer polynomials belong to the family of Jacobi polynomials. More specifically, we have
\begin{equation}\label{eq:Gegenbauer-Jacobi}
	\Gegenbauer{L}{\lambda}=\frac{\Gamma(L+2\lambda)\Gamma(\lambda+1/2)}{\Gamma(L+\lambda+1/2)\Gamma(2\lambda)}\Jacobi{L}{\lambda-1/2}{\lambda-1/2}.
\end{equation}
From \cite[Eq. (7.33.6)]{Szego1975} we have the following asymptotic estimates for the Gegenbauer polynomials:
\begin{equation}\label{eq:bounds-Gegenbauer}
	\Gegenbauer{L}{\lambda}=\begin{cases}
		\theta^{-\lambda}O(L^{\lambda-1}), & L^{-1}\leq \theta\leq \pi/2, \\
		O(L^{2\lambda-1}), & 0\leq \theta\leq L^{-1}.
	\end{cases}
\end{equation}
Near the endpoints, the asymptotic behavior of the Gegenbauer polynomials is given by the Mehler--Heine formula \eqref{eq:Mehler-Heine}. To find the exact expression for the polynomials $\Gegenbauer{L}{2}$, recall from \eqref{eq:Gegenbauer-Jacobi} that
\begin{equation}\label{eq:Gegenabuer-2-Jacobi}
	\Gegenbauer{L}{2}=\frac{\Gamma(L+4)\Gamma(5/2)}{\Gamma(L+5/2)\Gamma(4)}\Jacobi{L}{3/2}{3/2}\sim \frac{\sqrt{\pi}}{8}L^{3/2}\Jacobi{L}{3/2}{3/2}.
\end{equation}
Then, combining \eqref{eq:Mehler-Heine} and \eqref{eq:Gegenabuer-2-Jacobi}, we have
\begin{equation}\label{eq:Mehler-Heine-Gegenbauer-2}
	\lim_{L\tendsto\infty}L^{-3}\Gegenbauer{L}{2}\biggl(\cos\frac{z}{L}\biggr)=\frac{\sqrt{\pi}}{2\sqrt{2}}z^{-3/2}\BesselJ{3/2}(z).
\end{equation}

\subsubsection*{Hilb's formula for the Gegenbauer polynomials $\Gegenbauer{L}{2}$}

Using the relation \eqref{eq:Gegenabuer-2-Jacobi} and simplifying, Hilb's formula from \cref{thm:Hilb-Jacobi} takes the following form for the polynomials $\Gegenbauer{L}{2}$:
\begin{align*}
	\sin^{3/2}\theta\Gegenbauer{L}{2}(\cos\theta)=\sqrt{\frac{\pi}{8}}(L+3)(L+2)^{-1/2}(L+1)\frac{\theta^{1/2}}{\sin^{1/2}\theta}\BesselJ{3/2}\bigl((L+2)\theta\bigr)+\textup{Error},
\end{align*}
where
\begin{equation*}
	\textup{Error}=\begin{cases}
			\theta^{1/2}O(L^{-3/2}), & L^{-1}\leq \theta\leq \pi-\varepsilon,\\
			\theta^{7/2}O(L^{3/2}), & 0<\theta<L^{-1}.
		\end{cases}
\end{equation*}
Simplifying even more, we can write the previous formula as
\begin{equation}\label{eq:Hilb-Gegenbauer-2}
	\sin\theta\,\Gegenbauer{L}{2}(\cos\theta)=\sqrt{\frac{\pi}{8}}L^{3/2}(1+O(L^{-1}))\frac{\theta^{1/2}}{\sin\theta}\BesselJ{3/2}\bigl((L+2)\theta\bigr)+\textup{Error},
\end{equation}
where the error is
\begin{equation*}
	\textup{Error}=\begin{cases}
		O(L^{-3/2}), & L^{-1}\leq \theta\leq \pi-\varepsilon,\\
		\theta^{3}O(L^{3/2}), & 0<\theta<L^{-1}.
	\end{cases}
\end{equation*}

\subsection{Bessel functions}

From \cite[Section 3.8.5]{MagnusOberhettingerSoni1966}, we have the following identity:
\begin{equation}\label{eq:integral-Bessel}
	\int_{0}^{\infty}t^{-s-1}\BesselJ{\nu}(t)^2\dd t=\frac{1}{2^{s+1}}\frac{\Gamma(s+1)\Gamma(\nu-s/2)}{\Gamma(1+s/2)^2\Gamma(\nu+1+s/2)}.
\end{equation}
As a consequence, we can derive the following formula:
\begin{equation}\label{eq:integral-Bessel-log}
	\int_{0}^{\infty}\frac{\log t}{t}\BesselJ{\nu}(t)^2\dd t=\frac{1}{2\nu}\Bigl(\digamma(\nu)+\log 2\Bigr)+\frac{1}{4\nu^2}.
\end{equation}

%----------- BIBLIOGRAPHY

\end{document}